\documentclass[a4paper,11pt]{article}
\title{
From Cantor to Semi-hyperbolic Parameters \\
along External Rays\\
}
\author{
Yi-Chiuan Chen 
and
Tomoki Kawahira
}
\usepackage{enumerate}
\usepackage{graphicx}
\usepackage{amsmath}
\usepackage{amssymb}
\usepackage{amsthm}
\usepackage{color}
\theoremstyle{plain}
\newtheorem{thm}{Theorem}[section]
\newtheorem{prop}[thm]{Proposition}

\newtheorem{cor}[thm]{Corollary}
\theoremstyle{remark}
\newtheorem{remark}[thm]{Remark}
\theoremstyle{definition}
\newtheorem*{eg}{Example}
\newcommand{\C}{\mathbb{C}}
\newcommand{\Cbar}{\overline{\C}}
\newcommand{\R}{\mathbb{R}}
\newcommand{\D}{\mathbb{D}}
\newcommand{\Dbar}{\overline{\D}}
\newcommand{\Z}{\mathbb{Z}}
\newcommand{\N}{\mathbb{N}}
\newcommand{\M}{\mathbb{M}}
\newcommand{\T}{\mathbb{T}}
\newcommand{\X}{\mathbb{X}}
\newcommand{\abs}[1]{{\left| #1 \right|}}
\newcommand{\paren}[1]{{\left( #1 \right)}}
\newcommand{\brac}[1]{{\left\{ #1 \right\}}}
\newcommand{\cR}{{\mathcal{R}}}
\newcommand{\cV}{{\mathcal{V}}}
\newcommand{\st}{\,:\,}
\newcommand{\QED}{\hfill $\blacksquare$}

\newcommand{\e}{\epsilon}
\newcommand{\dist}{\mathrm{dist}\,}
\newcommand{\diam}{\mathrm{diam}\,}
\newcommand{\lam}{\lambda}
\newcommand{\Lam}{\Lambda}
\newcommand{\sS}{{\sf S}}
\newcommand{\sZ}{{\sf Z}}
\newcommand{\chat}{{\hat{c}}}

\begin{document}

\maketitle

\begin{abstract}
For the quadratic family $f_{c}(z) = z^2+c$ with $c$ in the exterior of the Mandelbrot set, 
it is known that every point in the Julia set moves holomorphically. 
Let $\hat{c}$ be a semi-hyperbolic parameter in the boundary of the Mandelbrot set. 
In this paper we prove that for each $z = z(c)$ in the Julia set, 
the derivative $dz(c)/dc$ is uniformly $O(1/\sqrt{|c-\hat{c}|})$ 
when $c$ belongs to a parameter ray that lands on $\hat{c}$.
We also characterize the degeneration of the dynamics along the parameter ray.
\end{abstract}

\section{Introduction and main results}
Let $\M$ be the {\it Mandelbrot set},
the connectedness locus of the quadratic family
$$
\brac{f_c:z \mapsto z^2 + c}_{c \, \in \, \C}.
$$ 
That is, the Julia set $J(f_c)$ is connected if and only if 
$c \in \M$.
For $c\not\in \M$, 
it is well-known that the Julia set $J(f_c)$ is a Cantor set, 
and the critical point $z=0$ does not belong to the Julia set.
Moreover, $f_c$ with $c \notin \M$ is {\it hyperbolic}:
i.e., there exist positive numbers $\alpha_c$ and $\beta_c$
such that $|Df_c^n(z)| \ge \alpha_c(1 + \beta_c)^n$ 
for any $n \ge 0$ and $z \in J(f_c)$.

\paragraph{Holomorphic motion of the Cantor Julia sets.}
For $c\not\in \M$, because of hyperbolicity, 
every point in $z \in J(f_c)$ moves holomorphically with $c$.
In other words, we have a {\it holomorphic motion} 
(\cite{BR, L, Mc, MSS}) of the Cantor Julia sets
over any simply connected domain in $\C-\M$. 
In this paper, we obtain some results regarding limiting behavior of this holomorphic motion 
when $c$ approaches the boundary of $\M$. 

Let us describe it more precisely:
For a technical reason, 
we consider the holomorphic motion of a Cantor Julia set 
over the topological disk $\mathbb{X}=\C - \M \cup \R_+$,
where $\R_{ + }$ denotes the set of positive real numbers.
For any base point $c_0 \in \X$,
there exists a unique map $H:\X \times J(f_{c_0})  \to \C$ such that
\begin{enumerate}[(1)]
\item
$H(c_0, z) = z$ for any $z \in J(f_{c_0})$;
\item
For any $ c \in \X$,
the map $z \mapsto {H}(c,z)$ 
is injective on $J(f_{c_0})$ and it extends to a quasiconformal map on $\Cbar$.
\item
For any $z_0 \in J(f_{c_0})$,
the map $c \mapsto H(c,{z_0})$ is holomorphic on $\X$.
\item
For any $c \in \X$,
the map $h_c(z):=H(c,z)$ satisfies $h_c(J(f_{c_0}))=J(f_c)$ and 
$f_c \circ h_c = h_c \circ f_{c_0}$ on $J(f_{c_0})$.
\end{enumerate}
See \cite[\S 4]{Mc} for more details.

\begin{figure}[htbp]
\begin{center}
\includegraphics[width=.60\textwidth]{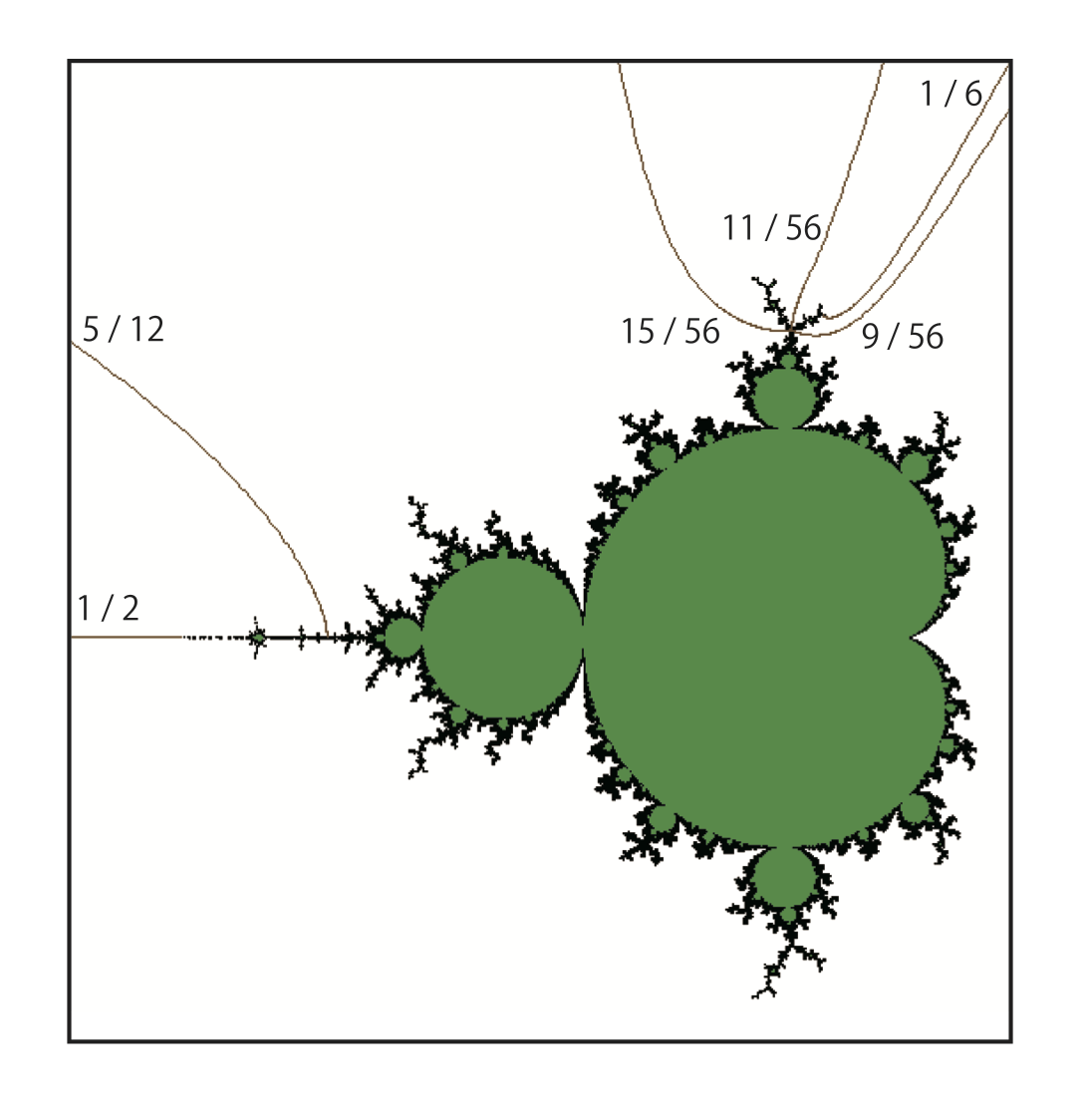}
\end{center}
\caption{The Mandelbrot set and the parameter rays of angles 
 $9/56$, $1/6$, $11/56$, $15/56$, $5/12$, and $1/2$.}
\label{fig_parameter_rays}
\end{figure}

\paragraph{Parameter rays.}
Let $\mathbb{D}$ denote the open disk of radius one centered at the origin.
 There is a unique biholomorphic function $\Phi_\M$ from 
$\Cbar- \M$ to $\Cbar- \Dbar$ 
satisfying $\Phi_\M(c)/c \to 1~(c \to \infty)$ with which the set
 \[ 
 \mathcal{R}_\M(\theta):=\{\Phi_\M^{-1}(r e^{i2\pi\theta}) \st 1< r < \infty\}
\]
 is defined and called the {\it parameter ray} of angle $\theta\in\mathbb{T}=\mathbb{R}/\mathbb{Z}$ of the Mandelbrot set $\M$. (This is a hyperbolic geodesic of the simply connected domain $\Cbar-\M$ starting at infinity.)  See Figure \ref{fig_parameter_rays}.
 Given $\theta$, 
 if the limit $\chat=\lim_{r\to 1^+}\Phi_\M^{-1}(re^{i2\pi\theta})$ 
 exists, then $\chat \in \partial \M$ is called the {\it landing point} of the parameter ray $\mathcal{R}_\M(\theta)$. We also say that $\theta$ is an {\it external angle} of the parameter $\chat$.

\begin{figure}[htbp]
\begin{center}
\includegraphics[width=.95\textwidth]{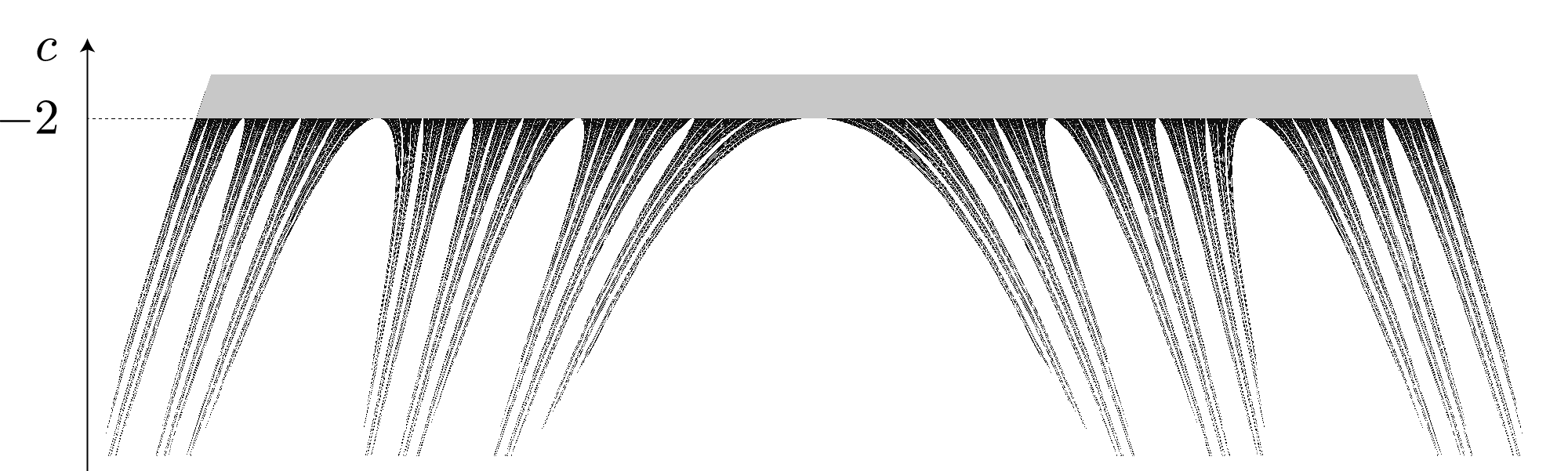}
\end{center}
\caption{
Each horizontal slice of the black part is the Julia set 
of parameter $c \in [-2.733, -2)$.
The gray part is the real slice of $J(f_c)$ for $c \in [-2, -1.875]$.
Note that $J(f_{-2})=[-2,2]$.
}\label{fig_realcantor}
\end{figure}

\begin{eg}[Real Cantor Julia sets]
When $c \notin \M$ approaches $\chat=-2$ along the real axis 
(equivalently, along the parameter ray of angle $1/2$),
$J(f_c)$ is contained in the real axis and its motion 
is depicted in Figure \ref{fig_realcantor}.
\end{eg}

\paragraph{Semi-hyperbolic parameters and Misiurewicz points.}
We are concerned with boundary behavior of the holomorphic motion given by the map $H$ above, along the parameter rays that land on a fairly large subset of $\partial \M$.

We say a parameter $\hat{c}$ in $\partial \M$ 
is {\it semi-hyperbolic} 
if the critical point is non-recurrent and belongs to the Julia set.
For each semi-hyperbolic parameter $\chat \in \partial \M$, 
there exists at least one parameter ray $\cR_\M(\theta)$ landing at $\chat$. 
(See \cite[Theorem 2]{D2}.
Indeed, there are at most finite number of parameter rays landing at $\chat$. 
See Remark \ref{rem:number_of_rays}.)
Note that for the quadratic polynomial $z^2+c$
(more generally, unicritical polynomials of the form $z^d+c$),
$\chat \in \partial \M$ being semi-hyperbolic implies  
it is a Collet-Eckmann parameter. 
(See \cite[Main Theorem \& p.51]{PRLS} also \cite[p.291 \& 299]{RL}.) 
Shishikura \cite{Shi} showed that for any open set $U$ intersecting with 
$\partial \M$, the semi-hyperbolic parameters in $U$ form a dense subset of Hausdorff dimension 2 of $U \cap \partial \M$.  
By a result of Douady \cite{D2}, 
the parameter ray $\cR_\M(\theta)$ 
lands on a semi-hyperbolic parameter if and only if 
$\theta \in \T$ is non-recurrent under the angle-doubling $t\mapsto 2t$ (mod $1$).
Hence every interval of $\T$ contains uncountably many angles for which
the parameter rays land on semi-hyperbolic parameters.
The geometric and dynamical properties of the Julia sets of 
semi-hyperbolic parameters 
are deeply investigated in a work of Carleson-Jones-Yoccoz 
\cite{CJY}. 
For example, if $\chat \in \partial \M$ is semi-hyperbolic,
then $J(f_\chat)$ is a locally connected dendrite 
such that $\Cbar-J(f_\chat)$ is a John domain.

\begin{figure}[htbp]
\begin{center}
\fbox{\includegraphics[width=.25\textwidth, bb = 0 0 600 600]{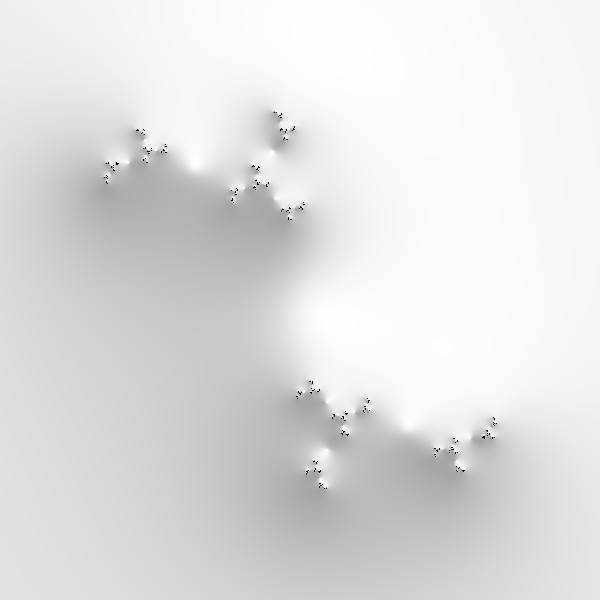}}
\fbox{\includegraphics[width=.25\textwidth, bb = 0 0 600 600]{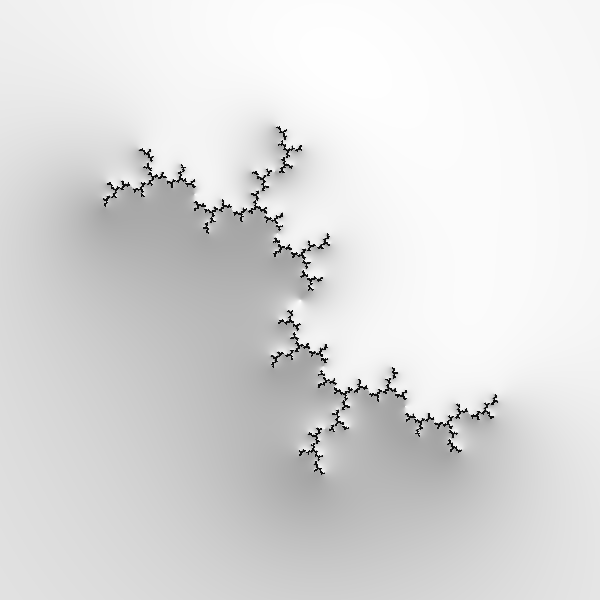}}
\fbox{\includegraphics[width=.25\textwidth, bb = 0 0 600 600]{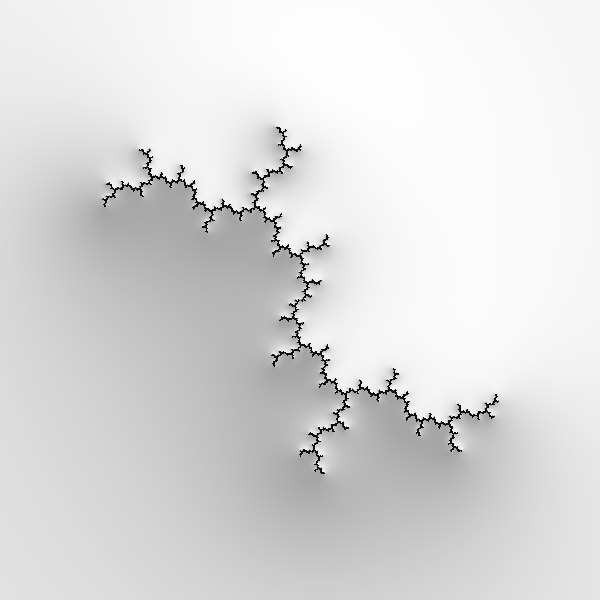}}
\\[.5em]
\fbox{\includegraphics[width=.25\textwidth, bb = 0 0 600 600]{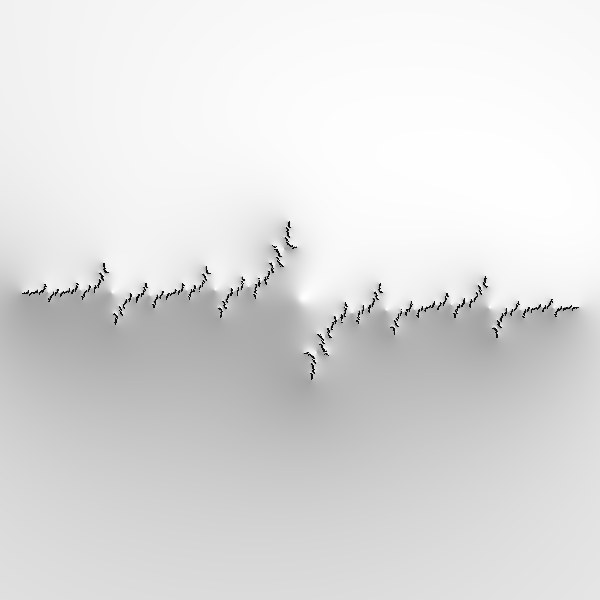}}
\fbox{\includegraphics[width=.25\textwidth, bb = 0 0 600 600]{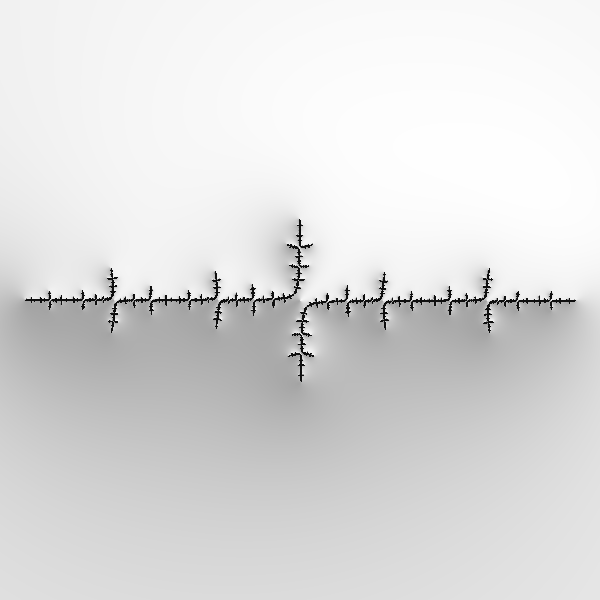}}
\fbox{\includegraphics[width=.25\textwidth, bb = 0 0 600 600]{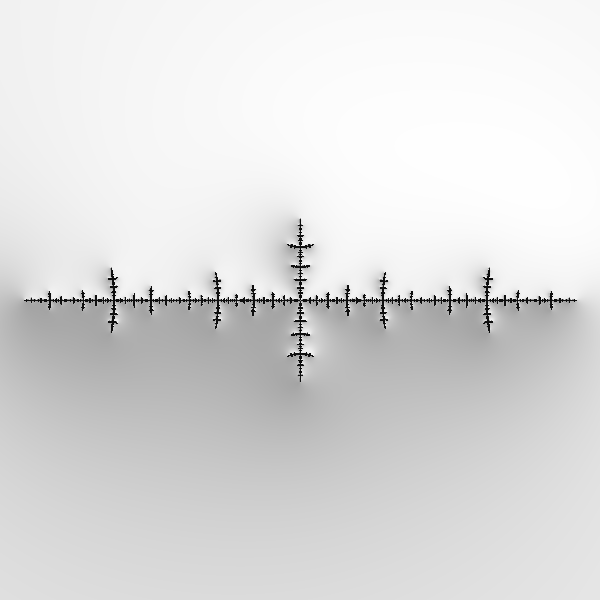}}
\\[.5em]
\fbox{\includegraphics[width=.25\textwidth,bb = 0 0 600 600]{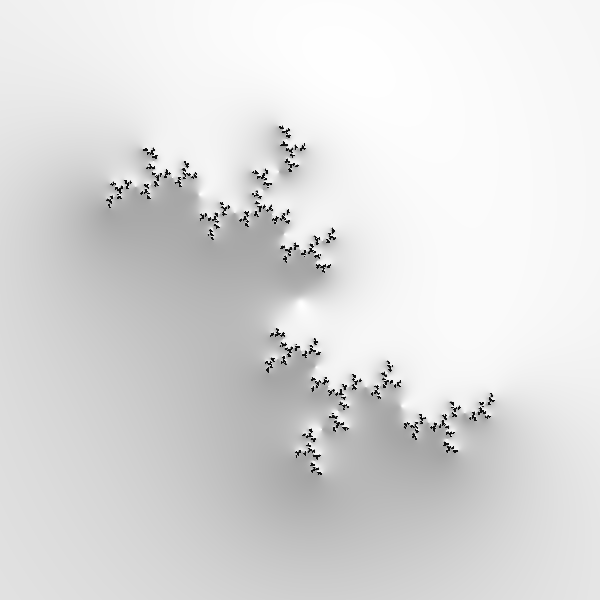}}
\fbox{\includegraphics[width=.25\textwidth,bb = 0 0 600 600]{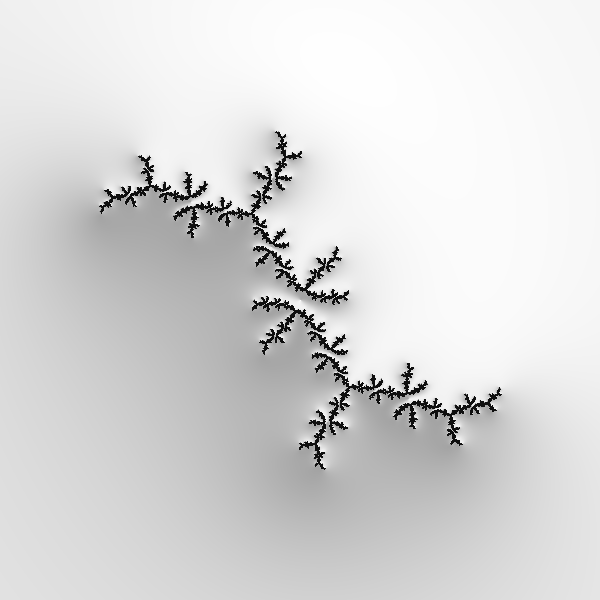}}
\fbox{\includegraphics[width=.25\textwidth,bb = 0 0 600 600]{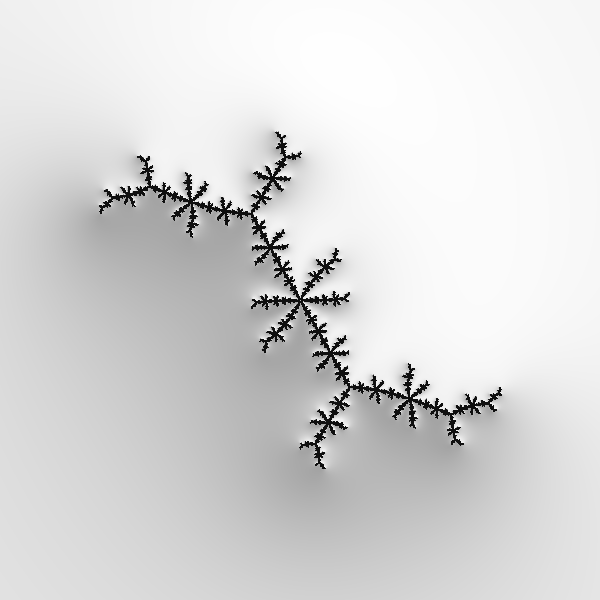}}
\\[.5em]
\fbox{\includegraphics[width=.25\textwidth,bb = 0 0 600 600]{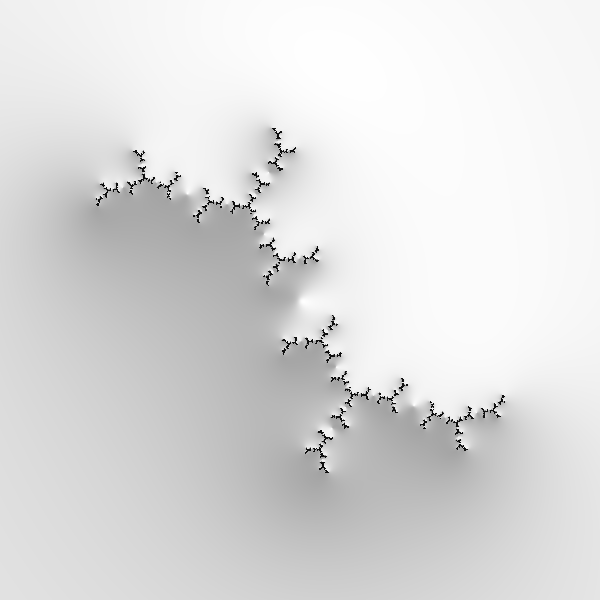}}
\fbox{\includegraphics[width=.25\textwidth,bb = 0 0 600 600]{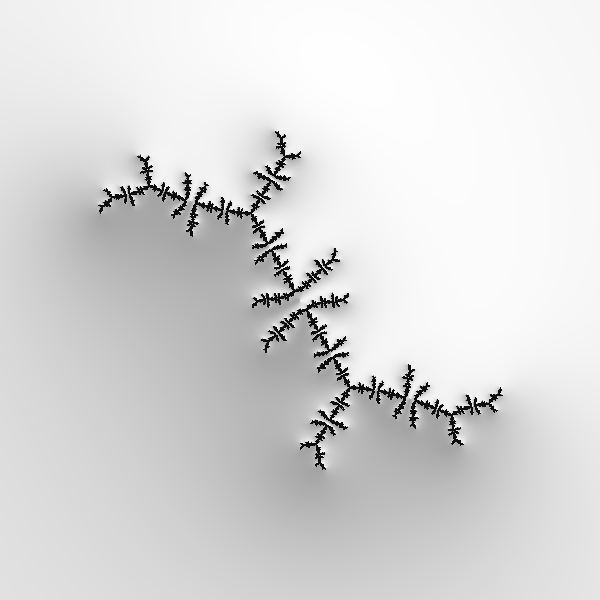}}
\fbox{\includegraphics[width=.25\textwidth,bb = 0 0 600 600]{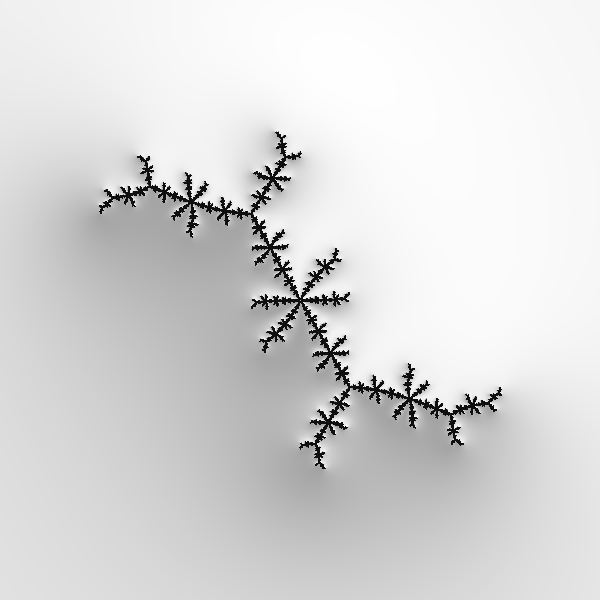}}
\\[.5em]
\fbox{\includegraphics[width=.25\textwidth,bb = 0 0 600 600]{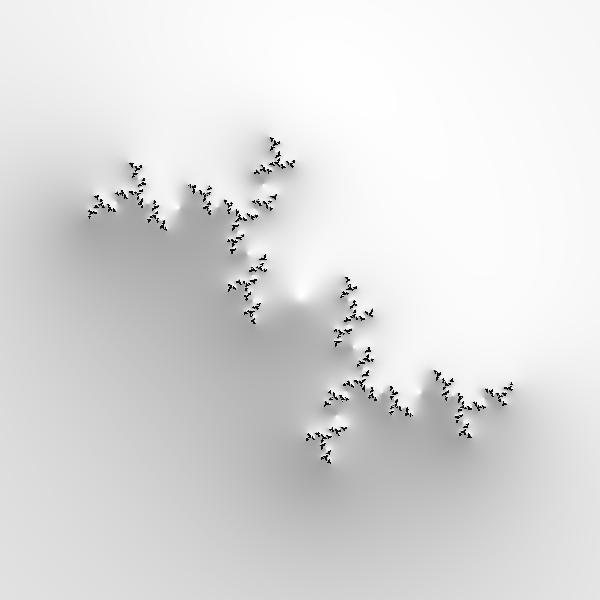}}
\fbox{\includegraphics[width=.25\textwidth,bb = 0 0 600 600]{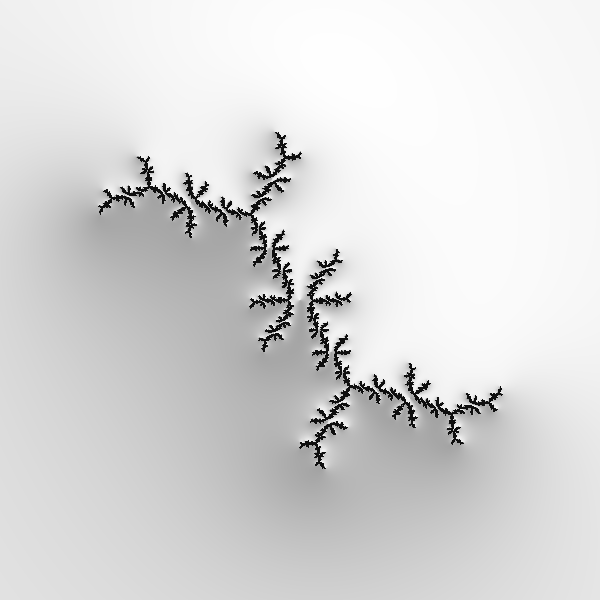}}
\fbox{\includegraphics[width=.25\textwidth,bb = 0 0 600 600]{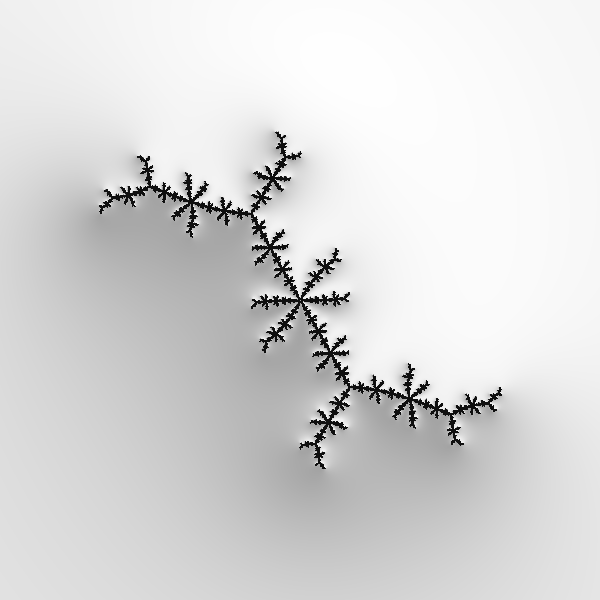}}
\end{center}
\caption{Holomorphic motion along the parameter rays of angles
$1/6$, $5/12$, $9/56$, $11/56$, and $15/56$.}
\label{fig_motions}
\end{figure}

A typical example of semi-hyperbolic parameter is 
a {Misiurewicz point}: 
We say a parameter $\chat$ is {\it Misiurewicz}
if the critical point of $f_{\hat{c}}$ is a pre-periodic point.
(By a {\it pre-periodic} point $z$ we mean 
$f_\chat^l(z)=f_\chat^{l+p}(z)$ for some integers $l$ and $p$ 
but $f_\chat^n(z) \not= z$ for all $n\ge 1$.)
It is known that such a Misiurewicz point $\hat{c}$ eventually lands on a repelling periodic cycle in the dynamics of $f_\chat$, and that the Misiurewicz points are contained in the boundary of the Mandelbrot set.
It is also known that the parameter $\hat{c}$ is Misiurewicz 
if and only if $\hat{c}$ is the landing point of 
$\mathcal{R}_\M(\theta)$ for some rational $\theta$ of even denominator. 
(See \cite[\'Expos\'e VIII]{DH1} and \cite[VIII, 6]{CG} for example.)
Holomorphic motions of the Julia sets along such rays are depicted in 
Figure \ref{fig_motions}.

\paragraph{Main results.}
Let $z_0$ be any point in $J(f_{c_0})$. 
Then the map $c \mapsto z(c):=H(c,z_0)$ is holomorphic over $\X=\C-\M \cup \R_+$. 
If we choose a semi-hyperbolic parameter $\chat \in \partial \M$, 
there exists a parameter ray $\cR_\M(\theta) \subset \X$ of 
angle $\theta \in \T-\brac{0}$ that lands on $\chat$.
As $c$ moves along the parameter ray $\cR_\M(\theta)$ toward $\chat$,
$z(c)=H(c,z_0)$ moves along an analytic curve in the plane. 

Our main theorem states that the speed of such a motion 
is uniformly bounded by a function of $|c-\chat|$: 

\begin{thm}[Main Theorem] \label{thm_main}
Let $\hat{c} \in \partial \M$ be a semi-hyperbolic parameter 
that is the landing point of $\mathcal{R}_\M(\theta)$.
Then there exists a constant $K>0$ that depends only on $\hat{c}$
such that for any $c\in\mathcal{R}_\M(\theta)$ 
sufficiently close to $\chat$ and any $z = z(c)\in J(f_c)$, 
the point $z(c)$ moves holomorphically with 
$$
\abs{\frac{dz(c)}{dc}}
\le 
\frac{K}{\sqrt{\abs{ c-\hat{c}}}}.
$$
\end{thm}

\medskip
The proof is given in Section \ref{sec:Proof of the main theorem}. 
By this theorem we obtain one-sided H\"older continuity of the holomorphic motion
along the parameter ray:

\begin{thm}[Holomorphic Motion Lands]\label{thm_HolomorphicMotionLands}
Let $\hat{c}\in \partial \M$ be a semi-hyperbolic parameter 
that is the landing point of $\mathcal{R}_\M(\theta)$, 
and let 
$c=c(r):=\Phi_\M^{-1}(re^{i2\pi\theta})$ with parameter $r \in (1,2]$.
Then for any $z(c(2))$ in $J(f_{c(2)})$, 
the improper integral 
$$
z(\hat{c}) := z(c(2))+ \lim_{\delta \to +0} \int_{2}^{1+\delta} \frac{dz(c)}{dc} ~\frac{dc(r)}{dr}~dr
$$
exists in the Julia set $J(f_{\chat})$. 
In particular, $z(c)$ is uniformly 
one-sided H\"older continuous of exponent $1/2$ at $c = \hat{c}$ along $\mathcal{R}_\M(\theta)$.
More precisely, there exists a constant $K'$ depending only on $\chat$ such that 
$$
|z(c)-z(\chat)| \le K' \sqrt{|c-\chat|}
$$
for any $c =c(r) \in \mathcal{R}_\M(\theta)$ with $1<r\le 2$.
\end{thm}
\medskip
This theorem implies:

\begin{thm}[From Cantor to Semi-hyperbolic]\label{thm_C2SH}
For any semi-hyperbolic parameter $\chat \in \partial \M$
and any parameter ray $\mathcal{R}_\M(\theta)$ 
landing at $\chat$, the conjugacy $H(c,\cdot)=h_c:J(f_{c_0}) \to J(f_c)$ 
converges uniformly to a semiconjugacy $h_{\chat}:J(f_{c_0}) \to J(f_\chat)$
from $f_{c_0}$ to $f_{\chat}$ as $c \to \chat$ along $\mathcal{R}_\M(\theta)$.
\end{thm}
The proofs of these theorems are given Section \ref{sec:Proofs of Theorems 1.2 and 1.3}. 
In Theorem \ref{thm_C2SH_2} below, 
we will specify where the semiconjugacy 
$h_{\chat}:J(f_{c_0}) \to J(f_\chat)$ fails to be injective.
Indeed, the semiconjugacy is injective except on a countable subset.

\medskip
By Theorems \ref{thm_HolomorphicMotionLands} and \ref{thm_C2SH},
we have a semiconjugacy $h_\chat \circ h_c^{-1}:J(f_c) \to J(f_\chat)$
with $|h_\chat \circ h_c^{-1}(z)-z| =O(\sqrt{|c-\chat|})$ as $c \to \chat$
along $\cR_\M(\theta)$. Thus we obtain:

\begin{cor}[Hausdorff Convergence]\label{cor_HausdorffConvergence}
The Hausdorff distance between $J(f_c)$ and $J(f_\chat)$ 
is $O(\sqrt{|c-\chat|})$ as $c \to \chat$ along $\cR_\M(\theta)$.
\end{cor}
This result is compatible with a result by Rivera-Letelier \cite{RL}. See Remark \ref{remark1}.

%%%%%%----------table of pictures
\fboxsep=0pt
\fboxrule=.5pt

\paragraph{Symbolic dynamics.}
Let 
$$
\Sigma_3:=\Big\{
{\bf s}=\{s_0,s_1,s_2,\ldots\} \st \ s_n=*, ~0\ \mbox{or}\ 1~ \mbox{for all}~ n\ge 0
\Big\}
$$ 
be the space consisting of sequences of $*$'s, $0$'s and $1$'s with the product topology, and  $\sigma$ be the left shift in $\Sigma_3$, 
$\sigma( {\bf s})={\bf s}'=(s'_0, s'_1, s'_2, \cdots)$ with $s'_i=s_{i+1}$. Let $$
\Sigma_2:=\Big\{
{\bf s}=\{s_0,s_1,s_2,\ldots\} \st \ s_n=0\ \mbox{or}\ 1~ \mbox{for all}~ n\ge 0
\Big\}\subset\Sigma_3
$$ 
be a closed subspace of $\Sigma_3$. 
A point ${\bf e}\in\Sigma_2$ is said to be {\it aperiodic} if $\sigma^n({\bf e})\not={\bf e}$ for any $n\ge 0$. 
  Two points ${\bf a}$ and ${\bf s}$ in $\Sigma_2$ are said to be {\it equivalent} with respect to aperiodic ${\bf e}\in \Sigma_2$, denoted by ${\bf a}\sim_{\bf e}{\bf s}$, if there is $k\ge 0$ such that $a_n=s_n$ for all $n\not= k$ and $\sigma^{k+1}({\bf a})=\sigma^{k+1}({\bf s})={\bf e}$.
It is plain to verify that the relation $\sim_{\bf e}$ is indeed an equivalence relation, and is the smallest equivalence relation that identifies $0{\bf e}$ with $1{\bf e}$.

Note that for $c \notin \M$ the dynamics of $f_c$ on the Julia set 
is conjugate to that of $\sigma$ on $\Sigma_2$.
We will use an aperiodic ${\bf e}$ to represent the (itinerary of the) 
non-recurrent critical orbit of the semi-hyperbolic $f_\chat$. 
Then ${\bf a}$ and ${\bf s}$ in $\Sigma_2$ are equivalent 
with respect to this ${\bf e}$
if and only if the points in $J(f_c)$ 
corresponding to ${\bf a}$ and ${\bf s}$ 
will degenerate to a point that eventually lands on 
the critical value $\chat$ in $J(f_\chat)$ as $c$ moves along 
the parameter ray landing on $\chat$.

Let $\mathcal{T}:\T\to\T$, $t\mapsto 2t ~(\bmod ~1)$ be the angle-doubling map. 
 Fix $\theta\in\T-\{0\}$, the two points $\theta/2$ and $(\theta+1)/2$ divide $\T$ into two open semi-circles $\T_0^\theta$ and $\T_1^\theta$ with $\theta\in\T_0^\theta$. Define the {\it itinerary} of a point $t$ under $\mathcal{T}$ with respect to $\theta$ as $\mathcal{E}^\theta (t)=\{ \mathcal{E}^\theta (t)_n \}_{n\ge 0}$
with
\[  \mathcal{E}^\theta (t)_n=
\begin{cases}
 0 & \mbox{for}~ \mathcal{T}^n(t)\in \T_0^\theta \\
 1 & \mbox{for}~ \mathcal{T}^n(t)\in \T_1^\theta \\
 * & \mbox{for}~ \mathcal{T}^n(t)\in \left\{\frac{\theta}{2},\frac{\theta+1}{2}\right\}.
\end{cases}
\]
The itinerary of $\theta$ itself, $\mathcal{E}^\theta (\theta)$, is called the {\it kneading sequence} of $\theta$.

Another consequence of Theorem \ref{thm_HolomorphicMotionLands} is as follows.
\begin{thm}[Symbolic Dynamics at Semi-hyperbolic Parameter] \label{corformain2}
Let $\chat$ be a semi-hyperbolic parameter with an external angle $\theta$ and ${\bf e}=\mathcal{E}^\theta (\theta)$ be the kneading sequence of $\theta$. Then $(J(f_\chat), f_\chat)$ is topologically conjugate to $(\Sigma_2/{\sim_{\bf e}}, \tilde\sigma)$, where $\tilde\sigma$ is induced by the shift transformation $\sigma$ of $\Sigma_3$.
\end{thm}

Theorem \ref{corformain2} also implies that the semiconjugacy in Theorem \ref{thm_C2SH} is one-to-one except at countable points where it is two-to-one.  

\begin{thm}[Almost Injectivity]\label{thm_C2SH_2}
Let $h_{\chat}:J(f_{c_0}) \to J(f_\chat)$ be the semiconjugacy given in Theorem \ref{thm_C2SH}.
For any $w \in J(f_\chat)$,
the preimage $h_\chat^{-1}(\{w\})$ has at most two points,
and it consists of two distinct points if and only if $f_\chat^n(w)=0$ for some $n \ge 0$.
\end{thm}
We prove these two theorems above in Section \ref{sec:Proofs of Theorems 1.5 and 1.6}.
More precise properties of the semiconjugacy can be found in Corollary \ref{corformain}.

\medskip

\paragraph{Structure of the paper.}
The structure of this paper is a little complicated,
but we belive this presentation requires less memory of the readers. 
In Section 2 we briefly summarize the notation and properties of 
the dynamics of $f_c(z)=z^2+c$ with semi-hyperbolic parameters. 
Section 3 is devoted for ``the derivative formula",
which is a key tool for our estimate.
In Section 4 we introduce the notion of ``Z-cycle" 
to describe the behavior of the orbits.
We also present Lemmas A, B, and C about Z-cycles,
whose proofs are given later.
In Section 5 we prove the Main Theorem by assuming these lemmas.
In Section 6 we introduce the notion of ``S-cycle"
and ``the S-cycle decompositions" of Z-cycles.  
We also present Lemmas A', B', and C',
whose proofs are given later as well.
Section 7 is devoted for Proposition S 
about stability of landing dynamic rays, 
and some lemmas that come from the assumption 
that the parameter $c$ moves along the parameter ray. 
In Section 8 we prove Theorems 1.2 and 1.3.
Then by assuming Lemmas A', B', and C', 
we prove Lemmas A and B in Sections 9 and 10 respectively.
Section 11 is devoted for some lemmas on hyperbolic metrics, 
and by using them, we prove Lemmas B', A', C', and C 
in Sections 12, 13, 14, and 15 respectively.
In Section 16 we work with symbolic dynamics, 
and finally in Section 17 we give proofs of Theorems 1.5 and 1.6.

\begin{remark} \label{remark1}
~~
\begin{itemize}
\item
The estimate in the Main Theorem is optimal. 
For example, if $\chat=-2$ (that is the Misiurewicz parameter with $f_\chat^2(0)=f_\chat^3(0)=2$),
then for $c=-2-\e$ with $\e >0$ 
the repelling fixed point on the positive real axis is given by 
$(1+\sqrt{9+4\e})/2=2+\e/3+o(\e)$.
Hence its preimages near the critical point are 
$z=\pm \sqrt{2\e/3}(1+o(\e))$, whose derivatives are 
$dz/d\e=\pm  (1/\sqrt{6\e})(1+o(\e))$. 
This implies that $|dz/dc|$ is compatible with $1/\sqrt{|c-\chat|}$.
See Figure \ref{fig_realcantor}.
\item
The results and the proofs in this paper are easily generalized to the unicritical family
$\brac{z \mapsto z^d+c \st c \in \C}$, simply by replacing the square root 
(``\,$\sqrt{|c-\chat|}$\,") by the $d$th root (``\,$|c-\chat|^{1/d}$\,")
in the Main Theorem. 
\item
In \cite{CK} the authors give a simple proof of the Main Theorem for $\chat=-2$.
\item
In \cite{D1}, Douady showed that the Julia set $J(f_c)$ 
continuously depends on $c$ at any semi-hyperbolic parameter $\chat$
in the sense of Hausdorff topology.
Moreover, in \cite{RL}, Rivera-Letelier showed that the Hausdorff distance 
between $J(f_c)$ and $J(f_{\chat})$ is $O(|c-\chat|^{1/2})$ for $c$ close enough to $\chat$, and that the Hausdorff dimension of the Julia set $J(f_c)$
converges to that of $\chat$ if $c$ tends to $\chat$ along the parameter ray. 
Our results, in addition, give the convergence of the dynamics.
\item
It is known that any parameter ray of odd denominator 
has a landing point $\chat$ on $\partial \M$ such that $f_\chat$
has a parabolic periodic point.
However, when $c$ moves along such a parameter ray,
$J(f_c)$ does not converge in the Hausdorff topology.
The discontinuity comes from the ``parabolic implosion",
which is also described in Douady's article \cite{D1}.
\item
Suppose $\chat\in\partial\M$ and $\chat\in J(f_\chat)$, and suppose $\chat$ has an external angle $\theta$. 
 There have been several results concerning the quotient dynamics for $f_\chat$ by kneading sequences.  If the kneading sequence $\mathcal{E}^\theta(\theta)$ is aperiodic, then the same statement as Theorem \ref{corformain2} that $(J(f_\chat), f_\chat)$ is topologically conjugate to $(\Sigma_2/{\sim_{\mathcal{E}^\theta(\theta)}}, \tilde\sigma)$ has been known by Bandt and Keller \cite{BK}.  Let $\approx_\theta$ be the smallest equivalence relation that if $t$, $t'$ are points in $\T$  such that for every $n$ either $\mathcal{E}^\theta(t)_n=\mathcal{E}^\theta(t')_n$ or $\mathcal{E}^\theta(t)_n=*$ or $\mathcal{E}^\theta(t')_n=*$, then $t$ is equivalent to $t'$. They also showed that $(J(f_\chat), f_\chat)$ is topologically conjugate to $(\T/{\approx_{\theta}}, \widetilde{\mathcal{T}})$ as well,
where $\widetilde{\mathcal{T}}$ is induced by the angle-doubling map  $\mathcal{T}$ on $\T/{\approx_{\theta}}$. 
Besides, for $f_\chat$  with locally connected  Julia set and  no irrational indifferent cycles, Kiwi \cite{K2} defined  $\equiv_\chat$ to be the smallest equivalence relation in $\T$ which identifies $t$ and $t'$ whenever the landing points of 
the dynamic rays $\mathcal{R}_\chat(t)$ and  $\mathcal{R}_\chat(t')$ coincide. 
(See Section \ref{sec:ParameterRayCondi} for the definition of the dynamic rays.) Then he showed that $(J(f_\chat), f_\chat)$ is topologically conjugate to $(\T/{\equiv_\chat}, \widehat{\mathcal{T}})$,
where $\widehat{\mathcal{T}}$ is induced by  $\mathcal{T}$ on $\T/{\equiv_\chat}$. (For $\chat$ a Misiurewicz parameter, Kiwi's result has been obtained earlier in \cite{AK}. However, in \cite{K2}  more general cases were considered including non-locally connected Julia sets.) 
\end{itemize}
\end{remark}

\section{Misiurewicz and semi-hyperbolic parameters}
In this section we briefly summarize the notation and properties of 
the dynamics of $f_c(z)=z^2+c$ with semi-hyperbolic parameters. 

\paragraph{Notation.}

\begin{itemize}
\item
Let $\N$ denote the set of positive integers. We denote the set of non-negative integers by $\N_0:=\{0\} \cup \N$.
\item
Let $\D(a,r)$ denote the disk in $\C$ centered at $a$ and of radius $r>0$.
When $a=0$ we denote it by $\D(r)$. 
\item
Let $\mathrm{N}(A, r)$ denote the open $r$-neighborhood of 
the set $A \subset \C$ for $r>0$. 
That is, $\mathrm{N}(A,r):=\bigcup_{a \in A} \D(a,r)$.
\item
For non-negative variables $X$ and $Y$, by $X \asymp Y$ we mean there exists 
an implicit constant $C>1$ independent of $X$ and $Y$ such that $X/C \le Y \le CX$.  
\item 
When we say ``for any $X \ll 1$" it means that 
``for any sufficiently small $X>0$". 
\item
Let $c$ be a parameter for the quadratic family 
$\{f_c(z)=z^2+c \st c \in \C\}$.
By $c \approx \chat$ we mean there exists an implicit constant $\delta>0$
independent of $c \neq \chat$ such that $|c-\chat|<\delta$. 
When we say ``the constant $C$ independent of $c \approx \chat$"
it means that $C$ does not depend on $c \neq \chat$ but it {\it may} 
depend on $\chat$.
\end{itemize}
 
\paragraph{Misiurewicz and semi-hyperbolic parameters.}
Let $\chat \in \partial \M$ be a Misiurewicz point 
with $f_\chat^l(0)=f_\chat^{l+p}(0)$ 
where we choose the  minimal $l$ and $p$ in $\N$. 
Then it is known that $f_\chat^l(0)$ is actually a repelling periodic point.

More generally, suppose that $\chat \in \partial \M$ is semi-hyperbolic,
and set $\hat{b}_n:=f_\chat^n(0)$ for each $n \ge 0$. 
Let $\Omega(\chat)$ denote the set of accumulation points of the set 
$\{\hat{b}_n\}_{n \ge 0}$,
i.e., the $\omega$-limit set of $0$. 
Moreover, by a result of Carleson, Jones, and Yoccoz \cite{CJY}, 
$\Omega(\chat)$ is a {\it hyperbolic set} in the sense of \cite{Shi}:
i.e., $\Omega(\chat)$ is compact; 
$f_\chat(\Omega(\chat)) \subset \Omega(\chat)$ (indeed, we have $f_\chat(\Omega(\chat))=\Omega(\chat)$);
and there exist constants $\alpha, \beta >0$ such that
$|Df_{\chat}^n(z)| \ge \alpha (1 + \beta)^n$ 
for any $z \in \Omega(\chat)$ and $n \ge 0$. 
For example, if $\chat$ is Misiurewicz, 
the set $\Omega(\chat)$ is the repelling cycle on which
the orbit of $0$ lands. 

For $\chat\in\partial\M$ a semi-hyperbolic parameter, it is proved in \cite{CJY} that there are constants $\epsilon>0$, $C>0$, and $0<\eta<1$ such that for all $z\in J(f_\chat)$, $n\ge 0$, and any connected component $B_n(z,\epsilon)$ of $f^{-n}_\chat(\D(z,\epsilon))$, we have 
\begin{equation}
 \diam B_n(z,\epsilon)<C \, \eta^n. \label{CJYcontraction}
\end{equation}  

In what follows we fix a $p \in \N$ such that $|Df_{\chat}^p(z)| \ge 3$
for any $z \in \Omega(\chat)$.
\footnote{Of course ``$3$" does not have particular meaning. Any constant bigger than one will do.}  
We first check:

\begin{prop}[Critical Orbit Lands]
\label{prop_crit_point_lands}
The critical orbit $\hat{b}_n =f_\chat^n(0)~(n \in \N_0)$ 
eventually lands on $\Omega(\chat)$. 
That is, 
there exists a minimal integer $l$ such that 
$\hat{b}_l =f_\chat^l(0) \in \Omega(\chat)$. 
\end{prop}

\paragraph{Proof.}
Suppose that $\hat{b}_n \notin \Omega(\chat)$ for every $n \in \N$.
Since $|Df_{\chat}^p(x)| \ge 3$ for any $x \in \Omega(\chat)$
we apply the Koebe distortion theorem (see \cite{Du}) 
to find a $\delta>0$ such that if  
$\hat{b}_n \in \mathrm{N}(\Omega(\chat), \delta)-\Omega(\chat)$, 
we have 
$$
\mathrm{dist}(\hat{b}_{n+p},\Omega(\chat)) \ge 2\, \mathrm{dist}(\hat{b}_n,\Omega(\chat)).
$$
(We also used compactness and invariance of $\Omega(\chat)$.
See also Remark \ref{remark2}.)
Hence there exists an accumulation point of the critical orbit
in $\overline{\C} - \mathrm{N}(\Omega(\chat), \delta)$.
However, it contradicts to the definition of $\Omega(\chat)$.
\QED

\medskip

Another remarkable fact is that the hyperbolic set $\Omega(\chat)$ 
moves holomorphically and preserves the dynamics (See \cite[\S 1]{Shi}):

\begin{prop}[Holomorphic Motion of $\Omega(\chat)$]
\label{prop_holo_motion_Omega}
There exist a neighborhood $\Delta$ of $\chat$ 
in the parameter plane $\C$
and a map $\chi:\Delta \times \Omega(\chat)  \to \C$
with the following properties:  
\begin{enumerate}[\rm (1)]
\item
$\chi(\chat,z) = z$ for any $z \in \Omega(\chat)$;
\item
For any $ c \in \Delta$,
the map $z \mapsto \chi(c,z)$ 
is injective on $\Omega(\chat)$ and it extends to a quasiconformal map on $\overline{\C}$.
\item
For any $z_0 \in \Omega(\chat)$,
the map $c \mapsto \chi(c,z_0)$ is holomorphic on $\Delta$.
\item
For any $c \in \Delta$,
the map $\chi_c(z):=\chi(c,z)$ satisfies 
$f_c \circ \chi_c = \chi_c \circ f_\chat$
on $\Omega(\chat)$.
\end{enumerate}
\end{prop}

\begin{figure}[htbp]
	\begin{center}
		\includegraphics[width=\textwidth, bb = 0 0 2329 1240]{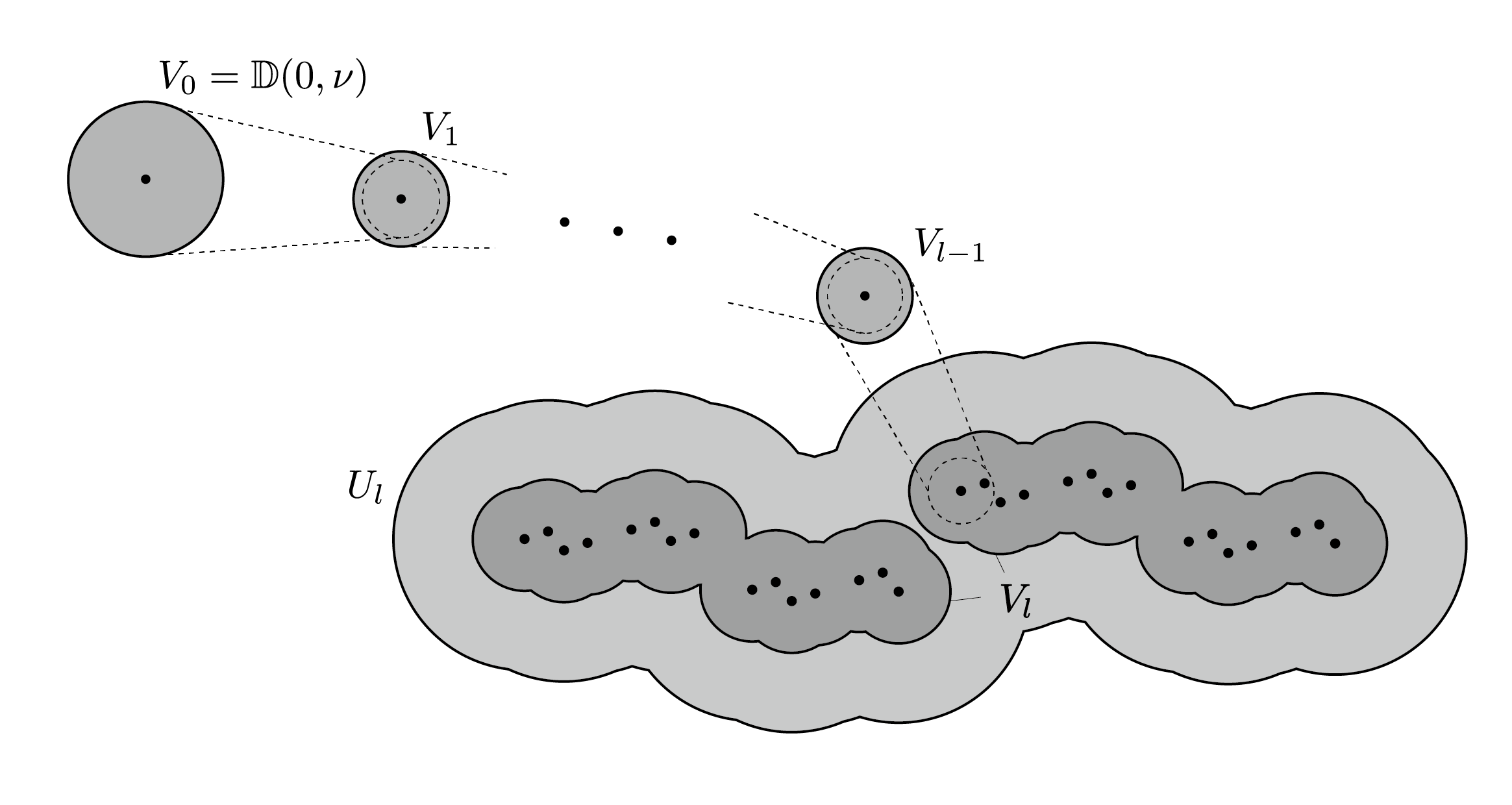}
	\end{center}
	\caption{$V_0$, $V_1$, $\cdots$, $V_l$ and $U_l$.}
	\label{Fig:V}
\end{figure}

\paragraph{Definition of $V_j$'s.}
Now we give a fundamental setting for the proofs of our results
that will be assumed in what follows.
\begin{itemize}
\item
Set $\Omega(c):=\chi_c(\Omega(\chat))$ 
for each $c \in \Delta$ given in Proposition \ref{prop_holo_motion_Omega}. 
Then $\Omega(c)$ is a hyperbolic subset of the Julia set $J(f_c)$. 
Since $J(f_c)$ is a Cantor set when $c \notin \M$, 
$\Omega(c)$ is a totally disconnected set for any $c \in \Delta$.
\item
Set $U_l:=\mathrm{N}(\Omega(\chat), R_l)$ for a sufficiently small $R_l>0$,
such that
\begin{itemize}
\item 
there is a constant $\mu \ge 2.5$ such that 
for any $c \approx \chat$ and $z \in U_l$ we have
$$
|Df_c^p(z)| \ge \mu;~~\text{and}
$$
\item
for any $c \approx \chat$, $U_l \Subset f_c^p(U_l)$.
\end{itemize}
Such an $R_l$ exists because $|Df_\chat^p(z)|\ge 3$ on $\Omega(\chat)$ 
and the function $(c,z) \mapsto |Df_c^p(z)|$ is continuous.
\item
We set $b_j(c):=\chi_c(\hat{b}_j) \in \Omega(c)$ for each $j \ge l$
and $c \in \Delta$.
By taking a smaller $\Delta$ if necessary,
we can also find an analytic family of 
pre-landing points 
$b_0(c)$, $b_1(c),\, \cdots, \,  b_{l-1}(c)$
over $\Delta$ 
such that $b_{j+1}(c)=f_c(b_j(c))$ and $\hat{b}_j=b_j(\chat)$ 
for each $j=0,1,\cdots, l-1$.
(For $j=0$, $b_0(c)$ is defined as a branch of $f_c^{-1}(b_1(c))$.)

\item
Choose disjoint topological disks 
$V_j$ for $j=0,1,\cdots,l-1$
such that
\begin{itemize}
\item
$V_0:=\D(0,\nu)$ for some $\nu \ll 1$.
We will add more conditions for $\nu$ later.
\item
For each $j=1,\cdots,l-1$, 
the topological disk $V_j$ contains 
$\hat{b}_j$ and satisfies $\diam V_j \asymp \nu^2$.
More precisely, there exists a constant $C_0>1$ 
independent of $j$ such that
$\nu^2/C_0 \le \diam V_j \le C_0\nu^2$.
\item
For any $c \approx \chat$ and each $j=0,1,\cdots,l-2$, 
we have $f_c(V_j) \Subset V_{j+1}$.
\end{itemize}
We also take a constant $C_0' >1$ such that for any $c \approx \chat$,
\begin{itemize}
\item
the set $V_l := \mathrm{N}(\Omega({\chat}), C_0' \nu^2)$ 
contains the topological disk $f_c(V_{l-1})$; and  
\item
at least for each $j=0, 1, \cdots, p-1$, 
$f_c^j(V_l) \Subset U_l$. 
\end{itemize}
We assume that $\nu$ is sufficiently small such that $V_j \cap V_l=\emptyset$ 
for each $j=0, 1, \cdots, l-1$.
Let $\cV$ denote the union $V_1 \cup V_2 \cup \cdots \cup  V_{l-1} \cup V_l$. 
See Figure \ref{Fig:V}.
\item
Let $\xi$ be the distance from $0$ to the closure of the set
$$
\{\hat{b}_1, \,\hat{b}_2,\, \cdots,\, \hat{b}_{l-1}\} \cup U_l.
$$
Since $0$ is not recurrent
(i.e., $0 \notin  \Omega(\chat)$),
we have $\xi>0$ if we take $R_l$ small enough.
We may assume in addition that $0 < \xi \le 1$ 
if we reset $\xi :=1$ when $\xi >1$.
If necessary, we replace $\nu$ 
so that $R_l$ and $C_0 \nu^2$ are smaller than $\xi/2$.
Then we have $|Df_c(z)|=2|z| \ge \xi$ 
for any $z \in \cV \cup U_l$ and $c \approx \chat$. 
\end{itemize}

\begin{remark} \label{remark2}
The backward dynamics of $f^p$ near $\Omega(\chat)$ is uniformly shrinking 
with respect to the Euclidean metric. 
For example, 
{\it one can find an $R>0$ depending only on $\chat$
such that for any $x \in \Omega(\chat)$
there exists a univalent branch $g$ of $f_\chat^{-p}$ 
on $\D(f_\chat^p(x),R)$ satisfying 
$g(f_\chat^p(x))=x$ and $g(\D(f_\chat^p(x),R)) \subset \D(x,R/2)$.}
Indeed, we first take an $R_0>0$ such that 
$f_\chat^p$ is univalent for any $\D(x,R_0)$ with $x \in \Omega(\chat)$.
By the Koebe $1/4$ theorem, 
$f_\chat^p(\D(x,R_0))$ contains $\D(f_\chat^p(x),R_0 |Df_\chat^p(x)| /4)$.
Since $|Df_\chat^p(x)| \ge 3$ on $\Omega(\chat)$,
there is a univalent branch $g$ of $f_\chat^{-p}$
on $\D(f_\chat^p(x), 3 R_0/4)$ with $g(f_\chat^p(x))=x$
and $|Dg(f_\chat^p(x))| \le 1/3$.
The Koebe distortion theorem implies that 
$g$ maps 
the disk $\D(f_\chat^p(x),R)$ 
into $\D(x,R/2)$ by taking a sufficiently small $R <3 R_0/4$.

We assume that the $R_l$ in the definition of $U_l$ 
is relatively smaller than this $R$, 
and we will implicitly 
apply this type of argument to 
the backward dynamics of $f_c$ near $U_l$ for $c \approx \chat$.
\end{remark}

\section{The derivative formula}

Recall that the map
$H:\mathbb{X} \times J(f_{c_0}) \to \C$
in Section 1
gives a holomorphic motion of the Julia set $J(f_{c_0})$
over the simply connected domain 
$\mathbb{X}=\C-\M \cup \R_+$ 
with the base point $c_0 \in \X$. 
For a given point $z_0 \in J(f_{c_0})$,
we want to have some estimates for the derivative 
of the holomorphic function $z(c)=H(c,z_0)$ at $c \in \X$.

In fact, such a holomorphic motion always exists for any
simply connected domain $\mathbb{Y}$ in $\C-\M$ with any base point $c_0 \in \mathbb{Y}$.
For a given $c \in \C-\M$, 
the derivative of such a motion at $c$
is independent of the choice of the domain 
$\mathbb{Y}$ containing $c$ and the basepoint $c_0$.
For example, it is convenient to consider the motion 
over the simply connected domain
$\mathbb{Y}:= \C -\M \cup \R_-$ (where $\R_-$ is the set of negative real numbers)
and assume that $\X$ and $\mathbb{Y}$ share the base point 
$c_0 \in \mathbb{Y} \cap \X = \C-\M \cup \R$.

Now we prove:

\begin{prop}
\label{prop_easy_estimate}
For any $c \notin \M$ and $z =z(c) \in J(f_c)$, 
we have
$$
\abs{\frac{d}{dc}{z}(c)} 
\le 
\frac{\, 1+\sqrt{1+6\,|c|\,}\,}{\dist(c, \partial \M)}.
$$
In particular, $\abs{dz/dc}=O(1/\sqrt{|c|})$ as $c \to \infty$.
\end{prop}

\paragraph{Proof.}
Let $\delta_c:=\dist(c, \partial \M)$ and $d_c:=(1+\sqrt{1+4\,|c|\,})/2$ for $c \in \C$.
Let $s(z):=\sup_{n \ge 0}~|f_c^n(z)|$ for $z \in J(f_c)$.
Since $f_c^{n+1}(z)=\paren{f_c^n(z)}^2+c$, we have $s(z) \ge s(z)^2-|c|$
and this implies $s(z) \le d_c$. 
Hence the Julia set $J(f_c)$ is contained in 
$\overline{\D(d_c)}$ for any $c \in \C$.

Now assume that $c \notin \M$. 
Then the disk $\D(c, \delta_c)$ is contained in either $\X=\C-\M\cup \R_+$ or $\mathbb{Y}=\C-\M\cup \R_-$,
and the motion of $J(f_{c_0})$ restricted to this disk is well-defined.
Let us consider a parameter $\zeta \in \D(c, \delta_c)$ such that $|\zeta-c| = \delta_c/2$.
Since $\delta_c \le |c|$, we have $|\zeta| \le 3|c|/2$ and thus 
the Julia set $J(f_\zeta)$
is contained in $\overline{\D(d_{3|c|/2})}$.
By applying the Cauchy integral formula, 
we obtain
\begin{align*}
\abs{
\frac{d}{dc}z(c)
}
=
\abs{
\frac{1}{2 \pi i} \int_{|\zeta-c| = \delta_c/2} 
\frac{z(\zeta)}{(\zeta-c)^2}\, d\zeta  
}
\le \frac{2 \,d_{3|c|/2}}{\delta_c}
=\frac{\, 1+\sqrt{1+6\,|c|\,}\,}{\dist(c, \partial \M)}.
\end{align*} 
Since $\M$ is contained in $\overline{\D(2)}$,
we have $|c|-2 \le \delta_c \le |c|$.
This implies $\abs{dz/dc}=O(1/\sqrt{|c|})$ as $c \to \infty$

\QED

\paragraph{The derivative formula.}
Our main theorem is based on the following formula (see also \cite{CKLY}):

\begin{prop}[The Derivative Formula]
\label{prop_Derivative_Formula}
For any $c \notin \M$ and $z =z(c) \in J(f_c)$, 
we have 
$$
\frac{d}{dc}{z}(c) = 
-\sum_{n = 1}^\infty\frac{1}{Df_c^{{n}}({z}(c))}.
$$
\end{prop}

\paragraph{Proof.}
Set $f:=f_c$ and ${z_n} = {z_n}(c): = f^{{n}}({z}(c))$. 
Then the relation ${z_{n + 1}} = {z_n}^2 + c$ implies
$$
\dfrac{d{z_{n + 1}}}{dc}
 = 2 {z_n}\cdot \dfrac{d{z_{n}}}{dc} + 1
\iff
\dfrac{d{z_{n}}}{dc}
 = 
 -\frac{1}{Df({z_n})} 
 + \frac{1}{Df({z_n})}  \dfrac{d{z_{n + 1}}}{dc}.
$$
Hence we have 
\begin{align*}
\frac{d}{dc}{z}(c) 
 = \dfrac{d{z_{0}}}{dc}
&=
 -\frac{1}{Df({z_0})} 
 + \frac{1}{Df({z_0})}  \dfrac{d{z_{1}}}{dc}\\
&= 
 -\frac{1}{Df({z_0})} + \frac{1}{Df({z_0})}  
\paren{-
 \frac{1}{Df({z_1})} + 
 \frac{1}{Df({z_1})}
   \dfrac{d{z_{2}}}{dc}}\\
&= 
 -\frac{1}{Df({z_0})}
  - \frac{1}{Df^{{2}}({z_0})}
 +  
\frac{1}{Df^{{2}}({z_0})}
\dfrac{d{z_{2}}}{dc}\\
&=
-\sum_{n = 1}^N\frac{1}{Df^{{n}}({z}(c))}
+\frac{1}{Df^{{N}}({z_0})}
\dfrac{d{z_{N}}}{dc}.
\end{align*}
By letting $N \to \infty$ we formally have the desired formula.
The series actually converges since 
$|dz_{N}/dc|$ is uniformly bounded 
by a constant depending only on $c$ (by Proposition \ref{prop_easy_estimate})
and $|Df^N(z_0)|$ grows exponentially by hyperbolicity of $f=f_c$. 
\QED

\medskip

\begin{remark}~
\begin{itemize}
\item
The estimate in Proposition \ref{prop_easy_estimate} 
is valid for any $c \in \C -\partial \M$.
Moreover, the derivative formula is also valid for any hyperbolic parameter in $\M$.
\item
Proposition \ref{prop_easy_estimate} implies an estimate
$$
\abs{\frac{dz}{dc}(c)}=O\paren{|c-\chat|^{-1-\beta}}
$$
if $c$ approaches $\chat \in \partial \M$ in such a way that 
$$
\mathrm{dist}(c,\partial \M) \ge C |c-\chat|^{1+\beta}
$$
for some constant $C>0$. 
The smallest possible value that $\beta$ can take is zero,
for example, when $c \to \chat=-2$ along the negative real axis.
Typically $\beta$ is positive, for example,
$\beta=1/2$ in the main theorem of \cite{RL}. 

\quad
In general, when $c$ approaches semi-hyperbolic 
$\chat \in \partial \M$ along a parameter ray landing at $\chat$,
it satisfies $\mathrm{dist}(c,\partial \M) \ge C |c-\chat|$
for some $C>0$, and thus $\beta=0$. 
(This is a combination of two facts: the John property 
of the complement of the Julia set $J(f_\chat)$ by \cite{CJY} and 
the asymptotic similarity between $J(f_\chat)$ and $\M$ 
at $\chat$ by \cite{RL}.)
This observation implies that our 
main theorem is stronger and it does not
come from the geometry of the Mandelbrot set.
We need the dynamics (the derivative formula) to prove it. 
\end{itemize}
\end{remark}

\section{Z-cycles}
For $c \approx \chat$, choose any $z=z_0 \in J(f_c)$.
The orbit $z_n:=f_c^n(z_0)~(n \in \N_0)$
may land on $V_0$ (or more precisely, on $V_0 \cap J(f_c)$),
and go out, then it may come back again. 
To describe the behavior of such an orbit, 
we introduce the notion of ``Z-cycle" for the orbit of $z$,
where ``Z" indicates that the orbit comes close to ``zero". 

We set $f:=f_c$ for brevity.

\paragraph{Definition {\normalfont (Z-cycle)}.}
A {\it finite Z-cycle} of the orbit $z_n=f^n(z_0)~(n \in \N_0)$ 
is a finite subset of $\N_0$ of the form
$$
\sZ=\brac{n \in \N_0 \st N \le n <N'}=[N,N') \cap \N_0, 
$$
such that 
 $z_N,\,z_{N'} \in V_0$ but 
$z_n \notin V_0$ if $N<n <N'$.
An {\it infinite Z-cycle} 
is an infinite subset of $\N_0$ of the form
$$
\sZ=\brac{n \in \N_0 \st N \le n <\infty}=[N,\infty) \cap \N_0, 
$$
such that 
 $z_N \in V_0$ but 
$z_n \notin V_0$ for all $n > N$.
By a {\it Z-cycle} we mean a finite or infinite Z-cycle.
In both cases, we denote them 
$\sZ=[N, N')$ or $\sZ=[N, \infty)$ 
for brevity.

\paragraph{Decomposition of the orbit by Z-cycles.}
For a given orbit $z_n=f^n(z_0)~(n \in \N_0)$ of $z_0 \in J(f_c)$,
the set $\N_0$ of indices is uniquely decomposed 
by using finite or infinite Z-cycles in one of the following three types:
\begin{itemize}
\item
The first type is of the form
\begin{equation}\label{eq_decomposition}
\N_0=[0,N_1) \sqcup \sZ_1 \sqcup \sZ_2 \sqcup \cdots,
\end{equation}
where $z_n \notin V_0$ for $n \in [0,N_1)$ and  
$\sZ_k:=[N_k, N_{k+1})$ is a finite Z-cycle for each $k \ge 1$.
\item
The second type is of the form
\begin{equation}\label{eq_decomposition2}
\N_0=[0,N_1) \sqcup \sZ_1 \sqcup \sZ_2 \sqcup \cdots \sqcup \sZ_{k_0},
\end{equation}
where $k_0 \ge 1$ such that $z_n \notin V_0$ for $n \in [0,N_1)$;  
$\sZ_k:=[N_k, N_{k+1})$ is a finite Z-cycle for each $1\le k \le k_0-1$;
and $\sZ_{k_0}=[N_{k_0}, \infty)$ is an infinite Z-cycle.
\item
The third type is just $\N_0=[0,N_1)$ with $N_1=\infty$,
where $z_n \notin V_0$ for all $n \in \N$. 
\end{itemize}
In the first and second types it is possible that $N_1=0$ and $[0,N_1)$ is empty.
For the second and third types, 
we set $\sZ_k:=\emptyset$ for any $k \ge 1$ for which $\sZ_k$ is not defined yet.
Hence we always assume that $\N_0$ 
formally has an infinite decomposition of the form (\ref{eq_decomposition})
associated with the orbit of $z_0 \in J(f_c)$.

\paragraph{The three lemmas.}
In what follows we assume the following ``parameter ray condition"
without (or with) mentioning:

\begin{quote}
{\bf ``Parameter ray condition"}.
{\it The parameter $c$ is always in the parameter ray $\cR_\M(\theta)$
that lands on $\chat$.
}
\end{quote} 
Now we present three principal lemmas about Z-cycle. 
(The proofs will be given later.) 

\paragraph{Lemma A.}
{\it 
There exists a constant $K_{\mathrm A}>0$
such that for any $c \approx \chat$, 
any $z=z_0 \in J(f_c)$, and for any 
Z-cycle $\sZ=[N,N')$ of the orbit $z_n=f_c^n(z)~(n \in \N_0)$, 
we have 
\begin{equation}
\sum_{i = 1}^{N'-N} \frac{1}{|Df_c^i(z_N)|}
\le 
\frac{K_{\mathrm A}}{\sqrt{|c-\chat|}},
\end{equation}
where we set $N'-N:=\infty$ if $N'=\infty$. 
}

\paragraph{Lemma B.}
{\it 
There exists a constant $K_{\mathrm B} > 0$
such that for any $c \approx \chat$
and any $N \le \infty$,
if $z=z_0 \in J(f_c)$ satisfies $z_n \notin V_0$ for any $n \in [0, N)$,
then we have \begin{equation}
\sum_{i = 1}^{N} \frac{1}{|Df_c^i(z_0)|}
\le K_{\mathrm B}.
\end{equation}
}
\medskip

In fact, $K_{\mathrm B}$ depends only on the choices of $\chat$ and $\nu$.
Hence we have:

\begin{cor}
\label{cor_of_Lemma B}
For any $c \approx \chat$ and any $z=z_0 \in J(f_c)$, 
if the orbit of $z$ 
never lands on $V_0=\D(\nu)$,
then the derivative satisfies
\begin{equation}
\abs{\frac{dz}{dc}}
\le 
\sum_{n = 1}^{\infty} 
\frac{1}{|Df_c^n(z_0)|}
\le K_{\mathrm B}.
\end{equation}
\end{cor}

\paragraph{Lemma C {\normalfont (Z-cycles Expand Uniformly)}.}
{\it 
There exists a constant $\Lam > 1$
such that for any $c \approx \chat$, 
any $z=z_0 \in J(f_c)$, and for any 
finite Z-cycle $\sZ=[N,N')$ of the orbit $z_n=f_c^n(z)~(n \in \N_0)$, 
we have 
\begin{equation}
|Df_c^{N'-N}(z_N)| \ge \Lam. 
\end{equation}
}

This $\Lam$ also depends only on the choice of $\nu$.
Indeed, $\Lam$ is bounded by a constant compatible with $\nu^{-1}$.

\section{Proof of the main theorem assuming Lemmas A, B, and C}
\label{sec:Proof of the main theorem}
We will use the derivative formula 
(Proposition \ref{prop_Derivative_Formula})
and Lemmas A, B, and C 
to show the inequality. 

For a given $c \approx \chat$
and $z=z_0 \in J(f_c)$, we consider the decomposition
$\N_0=[0, N_1)\sqcup \sZ_1\sqcup\sZ_2\sqcup\cdots$ 
as in (\ref{eq_decomposition}).
Set $f:=f_c$. Then we have
\begin{align*}
\abs{\frac{dz}{dc}}
&\le \sum_{n=1}^\infty \frac{1}{|Df^n(z_0)|}
=
\sum_{n=1}^{N_{1}}
\frac{1}{|Df^{n}(z_0)|}
+
\sum_{k \ge 1}\sum_{n \in \sZ_k}\frac{1}{|Df^{n+1}(z_0)|}\\
&=
\sum_{n=1}^{N_{1}}
\frac{1}{|Df^{n}(z_0)|}
+
\sum_{k \ge 1, \sZ_k \neq \emptyset}
\sum_{i=1}^{N_{k+1}-N_k}
\frac{1}{|Df^{N_k}(z_{0})|~|Df^i(z_{N_k})|}.
\end{align*}
By Lemma B, we obviously have $1/|Df^{N_1}(z_{0})| \le K_{\mathrm B}$.
By Lemma C, we have 
$$
|Df^{N_k}(z_{0})| =
|Df^{N_k-N_{k-1}}(z_{N_{k-1}})| 
\cdots
|Df^{N_2-N_{1}}(z_{N_{1}})| ~
|Df^{N_{1}}(z_{0})| 
\ge \Lam^{k-1}/K_{\mathrm B}
$$
as long as $\sZ_k \neq \emptyset$.
Hence by Lemma A, we have 
$$
\sum_{n=1}^\infty \frac{1}{|Df^n(z_0)|}
\le
K_{\mathrm B} + \sum_{k \ge 1} \frac{K_{\mathrm B}}{\Lam^{k-1}} 
\cdot \frac{K_{\mathrm A}}{\sqrt{|c-\chat|}}
=
K_{\mathrm B} + 
\frac{K_{\mathrm B}\Lam}{\Lam-1} 
\cdot \frac{K_{\mathrm A}}{\sqrt{|c-\chat|}}.
$$
We may assume that ${|c-\chat|} \le 1$ 
such that $K_{\mathrm B} \le K_{\mathrm B}/\sqrt{|c-\chat|}$.
Hence by setting 
$K:= K_{\mathrm B} + \dfrac{K_{\mathrm B}  K_{\mathrm A} \Lam}{\Lam -1}$, 
we have 
$
\abs{\dfrac{dz}{dc}} \le \dfrac{K}{\sqrt{|c-\chat|}}
$
for any $c \approx \chat$.

\QED

\section{S-cycles}
To show Lemmas A, B, and C, 
we introduce the notion of ``S-cycle". 

For $c \approx \chat$, set $f:=f_c$ and choose any $z=z_0 \in J(f_c)$.
The orbit $z_n:=f^n(z_0)~(n \in \N_0)$ may land on $\cV$.
Unless it lands exactly on the hyperbolic set $\Omega(c)$,
it will follow some orbit of $\Omega(c)$ for a while, 
and be repelled out of $U_l$ eventually.
Then it may come back to $\cV$, or land on $V_0$.
We define such a process as an ``S-cycle",
where ``S" indicates that orbit stays near the ``singularity" 
of the hyperbolic metric $\gamma$ to be defined in Section \ref{sec:Hyperbolic metrics}, or the cycle is relatively ``short" 
compared to Z-cycle. 

\paragraph{Definition {\normalfont (S-cycle)}.}
A {\it finite S-cycle} $\sS=[M,M')$ 
of the orbit $z_n=f^n(z_0)~(n \in \N_0)$ 
is a finite subset of $\N_0$ 
with the following properties:
\begin{itemize}
\item[(S1)]
$z_M \in V_j \subset \cV$
for some $j=1,2, \ldots, l$. If $M>0$ then $z_{M-1} \notin \cV$.
\item[(S2)]
There exists a minimal $m \ge 1$ such that 
for $n=M+(l-j)+mp$, 
$z_{n-p} \in U_l$ but $z_n \notin U_l$.
\item[(S3)]
$M'=M+(l-j)+mp+L$ for some $L \in [1,\infty)$
such that $z_n \notin V_0 \cup \cV$ 
for $n=M+(l-j)+mp+i~(0 \le i < L)$
and $z_{M'} \in V_0 \cup \cV$.
\end{itemize}
Note that in (S1), $z_{M-1}$ may be contained in $V_0$.
Note also that in (S2), some of $z_{n-p+1}, \, \cdots, z_{n-1}$ may {\it not} be contained in $U_l$.

An {\it infinite S-cycle} $\sS=[M,\infty)$ 
of the orbit $z_n=f^n(z_0)~(n \in \N_0)$ 
is an infinite subset of $\N_0$ satisfying either 
\begin{itemize}
\item 
Type {(I):}
(S1), (S2), and
\begin{itemize}
\item[(S3)']
$z_n \notin V_0 \cup \cV$ for all $n \ge M+(l-j)+mp$; 
\end{itemize}
\end{itemize}
or
\begin{itemize}
\item Type {(II):}
(S1) and
\begin{itemize}
\item[(S2)']
either $z_M =b_j(c)$ for $j<l$ or $z_M \in \Omega(c)$ for $j=l$. Equivalently, $z_n \in U_l$ 
for every $n= M + (l-j) +kp$ with $k \in \N$.
\end{itemize}
\end{itemize}

\paragraph{Decomposition of Z-cycles by S-cycles.}
Every Z-cycle $\sZ = [N,N')~(N \le \infty)$
of the orbit $z_n=f^n(z_0)~(n \in \N_0)$
has a unique decomposition by finite or infinite S-cycles.

For a finite Z-cycle $\sZ=[N,N')$,
there exists a finite decomposition 
$$
\sZ = \brac{N} \sqcup \sS_1 \sqcup \sS_2 \sqcup \cdots
\sqcup \sS_{k_0},
$$ 
where $\sS_k:=[M_k,M_{k+1})$ is a finite S-cycle 
for each $k=1,\, \cdots, \,k_0$ satisfying 
$N+1=M_1$ and $N'=M_{k_0+1}$. 

For an infinite Z-cycle $\sZ=[N,\infty)$,
there exists either a finite decomposition 
$$
\sZ = \brac{N} \sqcup \sS_1 \sqcup \sS_2 \sqcup \cdots
\sqcup \sS_{k_0},
$$
where $\sS_k:=[M_k,M_{k+1})$ is finite for 
$k=1,\, \cdots, \,k_0-1$ but infinite for $k=k_0$; or
an infinite decomposition 
$$
\sZ = \brac{N} \sqcup \sS_1 \sqcup \sS_2 \sqcup \cdots
$$
where $\sS_k:=[M_k,M_{k+1})$ is finite for any $k \ge 1$.

When we have a finite decomposition 
$
\sZ = \brac{N} \sqcup \sS_1 \sqcup \sS_2 \sqcup \cdots
\sqcup \sS_{k_0}
$,
we set $\sS_k:=\emptyset$ for $k > k_0$ and 
we assume that any Z-cycle formally has an infinite decomposition
 of the form 
$
\sZ = \brac{N} \sqcup \sS_1 \sqcup \sS_2 \sqcup \cdots
$. 
We call this {\it the S-cycle decomposition} of $\sZ$.

\paragraph{The three lemmas for S-cycles.}
Now we present three lemmas for S-cycles, that are parallel to 
Lemmas A, B, and C for Z-cycles:

\paragraph{Lemma A'.}
{\it 
There exists a constant $\kappa_{\mathrm A}>0$
such that for any $c \approx \chat$, 
any $z=z_0 \in J(f_c)$, and for any 
S-cycle $\sS=[M,M')$ of the orbit $z_n=f_c^n(z)~(n \in \N_0)$, 
we have 
\begin{equation}
\sum_{i = 1}^{M'-M} \frac{1}{|Df_c^i(z_M)|}
\le 
\kappa_{\mathrm A},
\end{equation}
where we set $M'-M:=\infty$ if $M'=\infty$. 
}

\paragraph{Lemma B'.}
{\it 
There exists a constant $\kappa_{\mathrm B} > 0$
such that for any $c \approx \chat$
and any $M \le \infty$,
if $z=z_0 \in J(f_c)$ satisfies $z_n \notin V_0 \cup \cV$ 
for $n \in [0, M)$,
then 
\begin{equation}
\sum_{i = 1}^{M} \frac{1}{|Df_c^i(z_0)|}
\le \kappa_{\mathrm B}.
\end{equation}
}

\paragraph{Lemma C' {\normalfont (S-cycles Expand Uniformly)}.}
{\it 
By choosing a sufficiently small $\nu$,
there exists a constant $\lam > 1$
such that for any $c \approx \chat$, 
any $z=z_0 \in J(f_c)$, and for any 
finite S-cycle $\sS=[M,M')$ 
of the orbit $z_n=f_c^n(z)~(n \in \N_0)$, 
we have 
\begin{equation}
|Df_c^{M'-M}(z_M)| \ge \lam. 
\end{equation}
}

\medskip
The proofs of these lemmas will be given later.

\section{Some lemmas concerning the parameter ray condition}\label{sec:ParameterRayCondi}
This section is devoted for some lemmas related to
the condition that $c$ is always on the parameter ray 
$\cR_\M(\theta)$ landing at $\chat$ (the ``parameter ray condition").

\paragraph{Dynamic rays for Cantor Julia sets. 
{\normalfont (See \cite[VIII, 3]{CG}, \cite[Appendix A]{M}.)}}
For any parameter $c \in \C$, the {\it B\"ottcher coordinate}
at infinity is a unique conformal map $\Phi_c$ 
defined near $\infty$ such that 
$\Phi_c (f_c(z))=\Phi_c(z)^2$ and $\Phi_c(z)/z \to 1$
as $z \to \infty$. 
Let $K(f_c)$ be the set of $z$ whose orbit 
is never captured in the domain of $\Phi_c$.
Then the boundary of $K(f_c)$ coincides with the Julia set $J(f_c)$.

 When $c \in \M$, the set $K(f_c)$ is connected and the B\"ottcher coordinate 
 extends to a conformal isomorphism 
 $\Phi_c:\C -K(f_c) \to \C -\overline{\D}$. 
 The {\it dynamic ray of angle} $t \in\T= \R/\Z$ 
is the analytic curve
$$
\cR_c(t) :=\brac{\Phi_c^{-1}(re^{2 \pi it}) \st r>1}.
$$
We say that {\it $\cR_c(t)$ lands at $z \in K(f_c)$} 
if $\Phi_c^{-1}(re^{2 \pi it})$ tends to $z$
as $r \searrow 1$.

When $c \notin \M$, the set $K(f_c)$ coincides with $J(f_c)$
which is a Cantor set.  
There exists a minimal $r_c>1$ such that the inverse $\Phi_c^{-1}$
extends to a conformal embedding of $\C-\overline{\D(r_c)}$ into $\C$
whose image contains the critical value $c=f_c(0)$.
(The Douady-Hubbard uniformization $\Phi_\M:\C-\M \to \C-\overline{\D}$ 
is given by setting $\Phi_\M(c):=\Phi_c(c)$.)
The dynamic ray of angle $t \in \T$ is partially defined in $\Phi_c^{-1}(\C-\overline{\D(r_c)})$,
and it extends to an analytic curve $\cR_c(t)$ 
landing at a point in $K(f_c)$ unless $2^n t=t_c$ for some $n \ge 1$,
where $t_c:=(2 \pi)^{-1} \arg \Phi_c(c)$.

\paragraph{Our setting and notation.}
Let us go back to our setting with semi-hyperbolic $\chat \in \partial \M$
where $\cR_\M(\theta)$ lands.
We will use the following facts and notations:
\begin{itemize}
\item
There is no interior point in $K(f_\chat)$ and thus $K(f_\chat)=J(f_\chat)$.
Moreover, $J(f_\chat)$ is connected and locally connected (\cite{CJY}). 
By Carath\'eodory's theorem, $\Phi_{\chat}^{-1}$ extends continuously
to $\C-\D$ and the dynamic ray $\cR_\chat(t)$ of any angle $t$ lands. 
\item
The angle $\theta$ is not recurrent under the angle doubling $t \mapsto 2 t$ (\cite[Thm.2]{D2}).
Set 
$$
\Theta:=\brac{2^{n+l-1} \theta \in \T \st n \ge 0}
$$
and let $\widehat{\Theta}$ denote its closure in $\T$, 
where $l$ is the minimal $l$ with $f^{l-1}_{\chat}(\chat)\in\Omega(\chat)$.
For $t \in \widehat{\Theta}$ the dynamic ray $\cR_\chat(t)$ lands on a point in 
the hyperbolic set $\Omega(\chat)$. 
(See Step 1 of Proposition S below.)
In particular, $\cR_\chat(2^{n+l-1} \theta)$ lands on $\hat{b}_{n+l} \in \Omega(\chat)$ 
for each $n \ge 0$.
\item
Let us fix an $r_0 >1$ and consider the compact set
$$
E_0:=\brac{
r e^{2 \pi i t} \st 
t \in \widehat{\Theta},\, r\in[r_0^{1/2^p},r_0]} \subset \C-\overline{\D}.
$$
By choosing $r_0$ close enough to $1$, 
the set $E(\chat):=\Phi^{-1}_\chat(E_0)$ is contained in $U_l$. 
\item
The parameter ray condition 
$c \in \cR_\M(\theta)$ is equivalent to $c \in \cR_c(\theta)$, 
or to $2 \pi t_c=\arg \Phi_c(c)= 2 \pi \theta$. 
Non-recurrence of $\theta$ assures that the dynamic rays
 $\cR_c(t)$ with $t \in \widehat{\Theta}$ 
are always defined and land on the Julia set.
\item
Since the B\"ottcher coordinate $\Phi_c(z)$ 
is holomorphic in both $c$ and $z$ as long as it is defined, 
$E(c):=\Phi^{-1}_c(E_0)$ is well-defined for each $c \approx \chat$ 
and also contained in $U_l$. 
More precisely, we choose the 
disk $\Delta$ in Proposition \ref{prop_holo_motion_Omega} small enough 
and assume that both $E(c)$ and $\Omega(c)$ moves holomorphically in $U_l$ for any $c \in \Delta.$  

\end{itemize}

Let us check the following proposition, that is interesting in its own right:

\paragraph{Proposition S {\normalfont (Stability of Landing Rays)}.}
{\it
For any $c \in \Delta$ (without assuming the parameter ray condition) and any $t \in \widehat{\Theta}$,
 the dynamic ray $\cR_c(t)$ lands on a point in the hyperbolic set $\Omega(c)$
 and $\cR_c(t) \cap U_l$ has uniformly bounded length.
In particular, $\cR_c(2^{n+l-1} \theta)$ lands on $b_{n+l}(c) \in \Omega(c)$ 
for each $n \ge 0$.
Moreover, the set
$$
\widehat{\cR}(c) 
:= \overline{\bigcup_{t \in \widehat{\Theta}} \cR_c(t)}
\subset \overline{\C}
$$
moves continuously in the Hausdorff topology on the Riemann sphere
as $c \to \chat$.
}

\paragraph{Proof.}
The proof breaks into three steps.

\paragraph{Step 1.}
We first consider the case of $c=\chat$. 
We claim: {\it
For any angle $t \in \widehat{\Theta}$,
the dynamic ray $\cR_\chat(t)$ 
lands on $\Omega(\chat)$ and $\cR_\chat(t)\cap U_l$ 
has uniformly bounded length. 
}

Let $x=x(t)$ denote the landing point of $\cR_\chat(t)$.
By the Carath\'eodory theorem,
$x(t)$ depends continuously on the angle $t$.
Since $x(2^{l-1} \theta)=\hat{b}_l \in \Omega(\chat)$
and any angle $t \in \widehat{\Theta}$ is an accumulation point of the orbit of $2^{l-1} \theta$
by the angle doubling, we obtain $x(t) \in \Omega(\chat)$.
(Note that $\Omega(\chat)$ is forward invariant and compact.)

Let us set $\cR:= \cR_\chat(t)$ and 
$$
\cR(n):=\brac{z \in \cR \st |\Phi_\chat(z)|^{2^{np}} \in [r_0^{1/2^p},r_0]}
$$
for $n \ge 0$ such that 
 $f_\chat^{np}(\cR(n))=f_\chat^{np}(\cR) \cap E(\chat)$ and the union
$$
\cR(0)\cup \cR(1) \cup \cR(2) \cup \cdots 
$$
coincides with the bounded 
arc $\cR-\Phi_\chat^{-1}(\brac{w \in \C \st |w|>r_0})$.
Note that the arc $f_\chat^{np}(\cR(n)) \subset E(\chat) \subset U_l$ 
has uniformly bounded length.
By the Koebe distortion theorem and the condition $|Df_\chat^{p}(z)| \ge \mu$ in $U_l$, 
we have 
$$
\mathrm{length} (\cR(n)) =O(\mu^{-n}),
$$
where the implicit constant is independent of the angle $t$.
Hence the dynamic ray $\cR$ has uniformly bounded length in $U_l$.

\paragraph{Step 2.}
Next we claim: {\it
For any $c \approx \chat$ and angle $t \in \widehat{\Theta}$,
the dynamic ray $\cR_c(t)$ 
lands on $\chi_c(x(t)) \in \Omega(c)$
and $\cR_c(t)\cap U_l$ 
has uniformly bounded length. 
}

Set $\cR':= \cR_c(t)$ and 
$$
\cR'(n):=\brac{z \in \cR' \st |\Phi_c(z)|^{2^{np}} \in [r_0^{1/2^p},r_0]}
$$
such that $f_c^{np}(\cR'(n)) =f_c^{np}(\cR') \cap E(c)$.
We also set $x':=\chi_c(x)$ where $x=x(t)$ is the landing point of 
$\cR=\cR_\chat(t)$ in $\Omega(\chat)$.
Since $\Omega(c)$ and $E(c)$ move holomorphically in $U_l$ 
with respect to $c \approx \chat$,
we may assume that the disk $D:=\D(f_\chat^{np}(x), R_l)$ contains
the point $f_c^{np}(x')=\chi_c(f_\chat^{np}(x))$
and the arcs $f_\chat^{np}(\cR(n))$ and $f_c^{np}(\cR'(n))$.
Since there exists a univalent branch $g_c$ of $f_c^{-np}$
defined on $D$ such that it sends $f_c^{np}(x')$ to $x'$
and $f_c^{np}(\cR'(n))$ to $\cR'(n)$, 
and since $|Df_c^{p}(z)| \ge \mu$ in $U_l$, we have 
$$
\dist(x', \cR'(n)) =O(\mu^{-n}).
$$
It follows that $\cR'=\cR_c(t)$ lands at $x'=\chi_c(x)$
and $\cR' \cap U_l$ has uniformly bounded length independent of 
$c \approx \chat$ and $t \in \widehat{\Theta}$. 

\paragraph{Step 3.}
Finally we show the continuity of the set $\widehat{\cR}(c)$.
It is enough to show:
{\it
For any $c \approx \chat$ there exists a 
homeomorphism $\phi_c: \widehat{\cR}(\chat) \to \widehat{\cR}(c)$
such that $\phi_c \to \mathrm{id}$ uniformly as $c \to \chat$
in the spherical metric.
}

By Step 2, the homeomorphism $\phi_c$ is naturally defined 
by $\phi_c(\infty)=\infty$, $\phi_c:=\chi_c$ on $\Omega(\chat)$, 
and $\phi_c:=\Phi_c^{-1} \circ \Phi_\chat$ on each ray $\cR_\chat(t)$ 
with $t \in \widehat{\Theta}$.

Now suppose that there exists an $\e>0$
such that for any $k \in \N$, we can find 
a pair of $c_k$ and $z_k$
such that $|c_k-\chat|\le 1/k$, $z_k \in \widehat{\cR}(\chat)$, 
and the spherical distance between 
$\phi_{c_k}(z_k)$ and $z_k$ exceeds $\e$. 
By taking a subsequence,
we may assume that $z_k$ has a limit $\zeta=\lim_{k \to \infty} z_k$ 
in $\widehat{\cR}(\chat)$.

Since the map $\Phi_c^{-1}(w)$ is continuous in both $c$ and $w$,
the map $\phi_c$ converges to identity as $c \to \chat$ 
locally uniformly
near each point of $\widehat{\cR}(\chat)-\Omega(\chat)\cup \{\infty\}$.
The convergence of $\phi_c$ near $\infty$ is uniform as well
in the spherical metric because $\Phi_c$ is tangent to identity near $\infty$.
Hence the limit $\zeta$ above belongs to $\Omega(\chat)$.

Let $W(n)$ denote the bounded subset of $\widehat{\cR}(\chat)$ 
given by
$$
W(n):= 
\Omega(\chat) \cup 
\bigcup_{t \in \widehat{\Theta}} 
\brac{\Phi_{\chat}^{-1}(r e^{2 \pi i t}) \st r \le r_0^{1/2^{np}}}.
$$
For any $n$, there exists a $k_n \in \N$ such that $z_k \in W(n)$
for any $k \ge k_n$. 
Now we define a point $x_k$ in $\Omega(\chat)$ as follow:
let $x_k:=z_k$ if $z_k \in \Omega(\chat)$.
Otherwise $z_k$ belongs to a dynamic ray 
$\cR_\chat(t_k)$ for some $t_k \in \widehat{\Theta}$,
and we let $x_k=x(t_k)$ be its landing point.
Then we obtain
\begin{align*}
|\phi_{c_k}(z_k)-z_k|
\le 
|\phi_{c_k}(z_k)-\phi_{c_k}(x_k)|
+
|\phi_{c_k}(x_k)-x_k|
+
|x_k-z_k|,
\end{align*}
where both $|\phi_{c_k}(z_k)-\phi_{c_k}(x_k)|$ and $|x_k-z_k|$ are $O(\mu^{-n})$ by Steps 1 and 2,
and $|\phi_{c_k}(x_k)-x_k|=|\chi_{c_k}(x_k)-x_k|=O(|c_k-\chat|)=O(1/k)$.
(See \cite[Corollary 2]{BR}.)
Hence $|\phi_{c_k}(z_k)-z_k|$ is bounded by $\e/2$ 
by taking sufficiently large $n$ and $k$. 
This is a contradiction.
\QED

\paragraph{}
The next lemma will be used in the proof of Lemma A:

\paragraph{Lemma T.}
{\it
Let $\hat{c}\in\partial \M$ be a semi-hyperbolic parameter.
There exists a positive constant
$C_{\mathrm T}=C_{\mathrm T}(\chat)$ 
such that $\dist (0,J(f_c))\ge C_{\mathrm T}\sqrt{|c-\chat|}$ for any $c \approx \chat$ on the 
parameter ray $\cR_\M({\theta})$ that lands at $\chat$. 
}

\paragraph{Proof.}
Since $f_c(z)-f_c(0)=(z-0)^2$,  it is equivalent to show
$$
\dist (c,J(f_c))\ge C_{\mathrm T}'{|c-\hat{c}|}
$$
for some constant $C_{\mathrm T}'=C_{\mathrm T}^2>0$ independent of $c \approx \chat$
with $c \in \cR_\M(\theta)$.

Set $a(c):=f_c^l(0)$ and $b(c):=b_l(c)$ for $c \approx \chat$.
Since $f_c^{l-1}$ is univalent near $c$, 
we have 
$$
\dist (c,J(f_c))\asymp \dist(a(c), J(f_c))
$$
by the Koebe distortion theorem.
By a result of Rivera-Letelier \cite[Appendix 2]{RL} and van Strien \cite[Theorem.1.1]{vS}
(see also Douady and Hubbard \cite[p.333, Lemma 1]{DH2} for Misiurewicz case), 
there exists a constant $B_0 \neq 0$
such that
$$
a(c)-b(c)=B_0(c-\chat) + O((c-\chat)^2)
$$ 
for $c \approx \chat$. 
Hence it is enough to show that there exists a constant $C_{\mathrm T}''>0$ such that 
\begin{equation}\label{eq_lem_Z}
\dist(a(c), J(f_c)) \ge C_{\mathrm T}'' |a(c)-b(c)|
\end{equation}
for $c \approx \chat$ with $c \in \cR_\M(\theta)$.

For each $z \in E(c)=\Phi_c^{-1}(E_0)$ defined in the proof of Proposition S, 
there exists an angle $t \in \widehat{\Theta}$ 
such that $\arg \Phi_c(z)=2 \pi t$. 
By Proposition S, the external ray $R_c(t)$ lands on a point $L_c(z)$
in $\Omega(c)$. 
Now we define a constant $\Gamma(c)$ for each $c \approx \chat$ by 
$$
\Gamma(c):=\inf \brac{\frac{\dist(z,J(f_c))}{|z-L_c(z)|} \in (0,1] \st z \in E(c) }
$$
and claim that its infimum 
$$
\Gamma:=\inf \brac{\Gamma(c) \st c \in \Delta}
$$
is a positive constant if we choose sufficiently small disk $\Delta$ centered at $\chat$.
Indeed, if there exists a sequence $c_k \to \chat$ 
in $\Delta$ such that $\Gamma(c_k) \to 0$,
then we have $\dist(z_k,J(f_{c_k})) \to 0$ for some $z_k \in E(c_k)$.
(Note that $|z-L_c(z)|$ is always bounded 
because $E(c)$ and $J(c)$ are uniformly bounded for $c \in \Delta$.)
However, it is impossible because $E(c)$ and $J(f_c)$ move continuously at $c =\chat$ and
$E(\chat)$ has a definite distance from $J(f_\chat)$.
Hence we obtain 
$$
\dist(z,J(f_c)) \ge \Gamma\, |z-L_c(z)|
$$
for each $z \in E(c)$ and $c \in  \Delta $.

Suppose that $c \in R_\M(\theta) \cap \Delta$
and $f_c^{np}(a(c)) \in E(c)$ for some $n \in \N_0$. 
Since $L_c(f_c^{np}(a(c))) = f_c^{np}(b(c))$, we have
$$
\dist(f_c^{np}(a(c)),J(f_c)) \ge 
\Gamma \, |f_c^{np}(a(c))-f_c^{np}(b(c))|.
$$
By Proposition S, if we choose sufficiently small $r_0$, 
then the length of the arc in the dynamic ray joining 
any $z \in E(c)$ and $L_c(z) \in \Omega(c)$ is uniformly and arbitrarily small. 
Thus there exists a univalent branch of $f_c^{-{np}}$ on the disk $\D(f_c^{np}(b(c)), 2 R_l)$
that sends both $f_c^{np}(a(c))$ and $f_c^{np}(b(c))$ to $a(c)$ and $b(c)$ respectively.
By the Koebe distortion theorem, we have (\ref{eq_lem_Z}).
\QED  

\begin{remark}
This proof is based on the argument to show that the basin at infinity of $f_\chat$ is a John domain. 
See \cite[\S 3]{CJY} and \cite[p.118]{CG}.
\end{remark}

\paragraph{}
The next lemma will be used in the proof of Lemma C:

\paragraph{Lemma U.}
{\it
There exists a constant $C_{\mathrm U} > 0$ with the following property:
for any $c \approx \chat$ with $c\in \cR_\M(\theta)$
and any $z_0 \in V_0 \cap J(f_c)$ such that $z_{n-p} \in U_l$ and 
$z_n \notin U_l$, we have $|Df_c^n(z_0)|\ge C_{\mathrm U}/|z_0|$.
}

\paragraph{Proof.}
By Lemma T (and its proof), 
we have $|z_0| \ge \dist(0, J(f_c)) \ge C_{\mathrm T} \sqrt{|c-\chat|}$
and $|b_0(c)| \asymp \sqrt{|b_l(c)-f_c^{l-1}(c)|} \asymp  \sqrt{|c-\chat|}$. 
Hence we have $|z_0| \ge C_1|b_0(c)|$ for some constant $C_1>0$
and it follows that 
$$
|z_1-b_1(c)|
=
|z_0^2-b_0(c)^2|  \le C_2 |z_0|^2
$$
where $C_2 :=1+C_1^2$. 

Now $z_n \notin U_l$ means that $|z_n-\hat{b}_n|\ge \dist (z_n, \Omega(\chat)) \ge R_l$.
Since $z_{n-p} \in U_l$, $z_n$ is still close to $\Omega(\chat)$
and by taking a smaller $R_l$ if necessary,
we may assume that there exists an $R>R_l$ 
independent of $c \approx \chat$ and 
$z_0 \in V_0 \cap J(f_c)$ such that $z_n \in \D(\hat{b}_n, R)$. 
Since we may assume that 
$|\hat{b}_n-b_n(c)| =|\hat{b}_n-\chi_c(\hat{b}_n)|\le R_l/2$ 
for any $c \approx \chat$, we have
$$
|z_n-b_n(c)| \ge |z_n-\hat{b}_n|-|\hat{b}_n-b_n(c)| \ge R_l/2.
$$
Let $G$ be a univalent branch of $f_c^{-(n-1)}$ 
defined on $\D(\hat{b}_n, 2R)$ 
(by taking smaller $R$ and $R_l$ if necessary)
that maps $b_n(c)$ to $b_1(c)$ 
and $z_n$ to $z_1$.
By the Koebe distortion theorem, we have
$$
|DG(z_n)| \asymp |DG(b_n(c))|
$$ 
and
$$ 
|z_1-b_1(c)|=|G(z_n)-G(b_n(c))| \asymp |DG(b_n(c))|~|z_n-b_n(c)|.  
$$
Since $|z_1-b_1(c)| \le C_2|z_0|^2$ and $|z_n-b_n(c)| \ge R_l/2$,
we have $|Df_c^{n-1}(z_1)| =|DG(z_n)|^{-1}\ge C_3/|z_0|^2$, where $C_3$ is a constant 
independent of $c \approx \chat$. 
Hence we have 
$$
|Df_c^n(z_0)|=|Df_c^{n-1}(z_1)|
~|Df_c(z_0)|\ge \frac{C_3}{|z_0|^2}\cdot (2 |z_0|) =\frac{2C_3}{|z_0|}.
$$
Set $C_{\mathrm U}:=2C_3.$ 
\QED\\

\paragraph{Geometry of the parameter ray.}

The following lemma will be used in the proof of Theorem \ref{thm_HolomorphicMotionLands}:

\paragraph{Lemma V.}
{\it
Let $\chat \in \partial \M$ be a semi-hyperbolic parameter
and $\cR_\M(\theta)$ a parameter ray landing on $\chat$.
Then the sequence $\brac{c_n}_{n \ge 0}$ in $\cR_\M(\theta)$ 
defined by 
$$
c_n:= \Phi_{\M}^{-1}\paren{r_0^{1/2^{np}} e^{2 \pi i \theta}}
$$
satisfies the following properties:
\begin{enumerate}[\rm (1)]
\item
$|c_{n + k}-\chat| =O(\mu^{-k}) |c_{n}-\chat|$ for any $n$ and $k \ge 0$.
\item
Let $\cR_\M(n)$ be the subarc of $\cR_\M(\theta)$ bounded by $c_n$ and $c_{n+1}$. Then 
$$
|c_{n+1}-c_n| \asymp \mathrm{length}( \cR_\M(n))  =O(\mu^{-n}).
$$
In particular, $\cR_\M(\theta)$ has finite length in a neighborhood of $\chat$.
\end{enumerate}
}

\paragraph{Proof.}
By a result by Rivera-Letelier \cite{RL}, 
there exists a constant $\hat{\lam} \neq 0$ such that
$
\Psi:=
\Phi_{\M}^{-1} \circ \Phi_{\chat}:\C-J(f_{\chat}) \to \C -\M
$ 
is of the form
$$
\Psi(z)
=\chat + \hat{\lam} (z-\chat)+O(|z-\chat|^{3/2})
$$
when $z \in \C-J(f_{\chat})$ and $z \approx \chat$.
In particular, $\Psi$ maps the dynamic ray $\cR_\chat(\theta)$ 
to the parameter ray $\cR_\M(\theta)$ conformally near the landing point $\chat$.
Hence it is enough to check that
{\it 
the points 
$$
z_n:= \Psi^{-1}(c_n)
=\Phi_{\chat}^{-1}\paren{r_0^{1/2^{np}} e^{2 \pi i \theta}}
$$
satisfies
\begin{enumerate}[\rm (1')]
\item
$|z_{n+k}-\chat| = O(\mu^{-k}) |z_{n}-\chat|$ for $k \ge 0$; and  
\item
the length of the subarc of $\cR_\chat(\theta)$ bounded by $z_n$ and $z_{n+1}$
is compatible with $|z_{n+1}-z_n|$ and is $O(\mu^{-n})$
\end{enumerate}
for sufficiently large $n$.
}

For each $t \in \widehat{\Theta}$ and $n \ge 0$, 
set 
$z_n(t):=\Phi_{\chat}^{-1}\paren{r_0^{1/2^{np}} e^{2 \pi i t}}$
such that the sequence $\{z_n(t)\}_{n \ge 0}$ converges along  
the external ray $\cR_\chat(t)$ to the landing point $x(t)$.
Note that $z_0(t)$ and $z_1(t)$ bound the arc 
$\cR_\chat(t) \cap E(\chat)$.
Since $E(\chat)$ and $\widehat{\Theta}$ are compact,
we have
\begin{itemize}
\item[(a)]
$|z_0(t)-x(t)| \asymp 1$; and 
\item[(b)]
$|z_0(t)-z_1(t)|  \asymp \mathrm{length} (\cR_\chat(t)  \cap E(\chat))$,
\end{itemize}
where the implicit constants are independent of $t \in \widehat{\Theta}$.

Now suppose that $n$ is large enough such that 
$np \ge l-1$ and thus $t_n:=2^{np}\theta \in \widehat{\Theta}$. 
Then we can find a univalent branch of $f_\chat^{-np}$
defined on a disk centered at $x(t_n)(=\hat{b}_{np+1})$
with a definite radius independent of $n$
that maps $z_0(t_n),\, z_k(t_n)$ and $x(t_n)$
univalently to $z_{n},\, z_{n+k}$ and $\chat$ respectively.
By the Koebe distortion theorem and (a) we have 
$$
\frac{|z_{n+k}-\chat|}{|z_n-\chat|}
\asymp
\frac{|z_k(t_n)-x(t_n)|}{|z_0(t_n)-x(t_n)|}
\asymp
|z_k(t_n)-x(t_n)|.
$$
We can find a univalent inverse branch $G_k$ of $f_\chat^{kp}$
defined on a disk centered at $x(t_{n+k})(=\hat{b}_{(n+k)p+1})$
with a definite radius independent of $n$ and $k$
that maps $z_0(t_{n+k})$ and $x(t_{n+k})$
univalently to $z_{k}(t_n)$ and $x(t_n)$.
Hence by Koebe again we have  
$$
|z_k(t_n)-x(t_n)| 
\asymp |DG_k(x(t_{n+k}))| \,|z_0(t_{n+k}) -x(t_{n+k})| =O(\mu^{-k}).
$$
It follows that $|z_{n+k}-\chat|= O(\mu^{-k}) |z_n-\chat|$
and we obtain (1'). 

By (b) and the same argument as above, 
the length of the subarc of $\cR_\chat(\theta)$ 
bounded by $z_n$ and $z_{n+1}$ is 
uniformly compatible with $|z_{n+1}-z_n|$ for any $n \ge 0$.
As a corollary of Step 1 of Proposition S,
we conclude that the length is $O(\mu^{-n})$. Thus we obtain (2').
\QED

\medskip

\begin{remark}\label{rem:number_of_rays}
Since there exist at most finitely many dynamic 
rays of the Julia set $J(f_\chat)$ landing at $\chat$ 
(see Thurston \cite[Theorem II.5.2]{Th} or Kiwi \cite[Theorem 1.1]{K1}),
the asymptotic similarity between $J(f_\chat)$ and 
$\M$ at $\chat$ by Rivera-Letelier \cite{RL}
implies that $\M$ has the same finite number of parameter rays
landing at $\chat$. (cf. \cite[VIII, 6]{CG}. See also \cite[Chapter 6]{Mc}.)
\end{remark}

\section{Proofs of Theorems \ref{thm_HolomorphicMotionLands} and \ref{thm_C2SH}}
\label{sec:Proofs of Theorems 1.2 and 1.3}

\paragraph{Proof of Theorem \ref{thm_HolomorphicMotionLands}.}
We combine the Main Theorem and Lemma V.
It is enough to show the existence of the improper integral
$$
z(c(r_0))+\lim_{\delta \to +0} \int_{r_0}^{1+\delta} \frac{dz(c)}{dc}\frac{dc(r)}{dr} \, dr
=
z(c(r_0))+\sum_{n\ge 0} 
\int_{\cR_\M(n)}
\frac{dz(c)}{dc} \, dc,
$$
where $r_0>1$ is a constant given in the definition of the set $E_0$
in the previous section, 
and $\cR_\M(n)$ is the subarc of $\cR_\M(\theta)$ bounded by $c_{n}$ and $c_{n+1}$ 
defined in Lemma V.
Note that by Lemma V, we obtain
$$
\mathrm{length}\cR_\M(n)  \asymp 
|c_{n+1}-c_n| \le |c_{n+1}-\chat|+|c_{n}-\chat|=O(|c_{n}-\chat|)
$$
and
$$
|c_{n}-\chat| \le \sum_{m \ge n}\mathrm{length}\cR_\M(m) =O(\mu^{-n}).
$$
Note also that 
\begin{equation}\label{eq_c_n}
|c_n -\chat| \asymp |c-\chat|
\end{equation}
for any $c \in \cR_\M(n)$, where the implicit constant is independent of $n$
by the Koebe distortion theorem, 
applied in the same way as the proof of Lemma V.

By the Main Theorem we obtain
\begin{align*}
\sum_{n\ge 0} 
\int_{\cR_\M(n)}
\abs{\frac{dz(c)}{dc}} \, |dc|
\le&
\sum_{n\ge 0} 
\int_{\cR_\M(n)}
\frac{K}{\sqrt{|c-\chat|}} \, |dc|\\
\asymp&
\sum_{n\ge 0} 
\frac{K}{\sqrt{|c_n-\chat|}} \, \mathrm{length}\cR_\M(n)\\
=& 
\sum_{n\ge 0}
O\paren{ 
\frac{1}{\sqrt{|c_n-\chat|}} \, |c_{n}-\chat|
}\\
=&
\sum_{n\ge 0} O(\mu^{-n/2})<\infty.
\end{align*}
Hence the improper integral above converges absolutely to some $z(\chat)$. 

To show the one-sided H\"older continuity,
it is enough to check
$|z(c_n) -z(\chat)| =O(\sqrt{|c_n-\chat|})$
for each $c_n$ by \eqref{eq_c_n}.
The same argument as above yields
$$
|z(c_n) -z(\chat)|
\le 
\sum_{k \ge 0} 
\int_{\cR_\M(n+k)}
\abs{\frac{dz(c)}{dc}} \, |dc|
\le
\sum_{k\ge 0} O(\sqrt{|c_{n+k}-\chat|}).
$$
By (1) of Lemma V, we have 
$|c_{n+k}-\chat|=O(\mu^{-k})|c_{n}-\chat|$
for each $k \ge 0$ 
and thus $|z(c_n) -z(\chat)|=\sum_{k\ge 0} O(\mu^{-k/2})\sqrt{|c_{n}-\chat|}
=O(\sqrt{|c_{n}-\chat|})$.

Since it is clear that  $z(\hat{c})$ is confined in a bounded region, to show $z(\hat{c}) \in J(f_{\hat{c}})$, we only need to show $\lim_{c\to\hat{c}}(z(c)^2+c)=(\lim_{c\to\hat{c}}z(c))^2+\lim_{c\to\hat{c}}c$, but this follows from the continuity of the quadratic map. 
\QED\\

\paragraph{Proof of Theorem \ref{thm_C2SH}}
For each $z_0 \in J(f_{c_0})$ and 
its motion $z(c)=h_c(z_0)=H(c,z_0)$ along 
the parameter ray $\cR_\M(\theta)$, 
we define $h_\chat(z_0)$ by the limit $z(\chat)$ given 
in Theorem \ref{thm_HolomorphicMotionLands}.
Since $h_c$ is continuous and the convergence of $h_c$ to $h_\chat$ 
as $c \to \chat$ along the parameter ray $\mathcal{R}_\M(\theta)$ 
is uniform, $h_\chat$ is continuous as well.
Hence $f_\chat \circ h_\chat = h_\chat \circ f_{c_0}$ is obvious and 
it is enough to show the surjectivity of $h_{\chat}:J(f_{c_0}) \to J(f_\chat)$.
First we take any repelling periodic point $x \in J(f_\chat)$.
Since there is a holomorphic family $x(c)$ of repelling periodic points 
for $c$ sufficiently close to $\chat$ such that $x=x(\chat)$, 
we have some $z_0 \in J(f_{c_0})$ with 
$h_c(z_0)=x(c)$ for any $c \approx \chat$ with $c \in \cR_\M(\theta)$. 
In particular, we have $h_\chat(z_0)=x.$
Next we take any $w \in J(f_{\chat})$ and a sequence of repelling 
periodic points $x_n$ of $f_{\chat}$ that converges to $w$ as $n \to \infty$.
(Such a sequence exists since repelling periodic points are dense in the Julia set.)
Let $z_n \in J(f_{c_0})$ be the repelling periodic point with
$h_\chat(z_n)=x_n$. 
Then any accumulation point $y$ of the sequence $z_n$ 
satisfies $h_\chat(y)=w$ by continuity.\QED

\section{Proof of Lemma A assuming Lemmas A' and C'}

Without loss of generality we may assume that $N=0$,
i.e., $z=z_0 \in V_0 \cap J(f_c)$.
We set $f:=f_c$. 
Now consider the S-cycle decomposition 
$
\sZ = \brac{0} \sqcup \sS_1 \sqcup \sS_2 \sqcup \cdots
$
of $\sZ=[0,N')$ where $\sS_k=[M_k,M_{k+1})$ if $\sS_k \neq \emptyset$,
 and $M_1=1$. 
Then we have 
\begin{align*}
\sum_{i=1}^{N'}\frac{1}{|Df^i(z)|}
&= \frac{1}{|Df(z)|}+
\sum_{k\ge 1}\sum_{n \in \sS_k}\frac{1}{|Df^{n+1}(z)|}\\
&= \frac{1}{2|z|}+
\sum_{k \ge 1,\, \sS_k \neq \emptyset}
\sum_{i=1}^{M_{k+1}-M_k}\frac{1}{|Df^i(z_{M_k})|~|Df^{M_k}(z)|}\\
&\le \frac{1}{2|z|}+
\sum_{k \ge 1,\, \sS_k \neq \emptyset}
\frac{\kappa_{\mathrm A}}{|Df^{M_k}(z)|}
\end{align*}
by Lemma A'.
If $\sS_k \neq \emptyset$, then by Lemma C',
$$
|Df^{M_k}(z)|
=
|Df^{M_k-M_{k-1}}(z_{M_{k-1}})|
\cdots
|Df^{M_2-M_{1}}(z_{M_1})|~|Df(z)|
\ge \lam^{k-1}\cdot 2|z|,
$$
where $M_1=1$. 
Hence we have $|Df^{M_k}(z)|^{-1} \le 1/(\lam^{k-1}\cdot 2|z|)$ for any $k$.
Moreover, by Lemma T, 
we have $\dist(0,J(f_c)) \ge C_{\mathrm T} \sqrt{|c-\chat|}$ 
for $c \approx \chat$ on the parameter ray, and thus
\begin{align*}
\sum_{i=1}^{N'}\frac{1}{|Df^i(z)|}
&\le
\frac{1}{2|z|}+
\sum_{k=1}^\infty \frac{\kappa_{\mathrm A}}{\lam^{k-1}\cdot (2|z|)}\\
&\le
\frac{1}{2\cdot\dist(0, J(f_c))}
\brac{1+
\sum_{k=1}^\infty \frac{\kappa_{\mathrm A}}{\lam^{k-1}}
}\\
&\le
\frac{1}{2C_{\mathrm T}\sqrt{|c-\chat|}}
\brac{1+\kappa_{\mathrm A}\frac{\lam}{\lam-1}}.
\end{align*}
Hence by setting $K_{\mathrm A}:=(2C_{\mathrm T})^{-1}\brac{1+ \kappa_{\mathrm A}\lam /(\lam -1)}$, we have the claim.
\QED

\section{Proof of Lemma B assuming Lemmas A', B' and C'}
Just like the S-cycle decompositions of Z-cycles, 
we have a finite or infinite decomposition of the form
$$
[0,N) = [0,M_1) \sqcup \sS_1 \sqcup \sS_2 \sqcup \cdots 
$$
where we have the following three cases:
\begin{enumerate}
\item
$N=M_1 \le \infty$ and $z_n \notin V_0 \cup\cV$ for any $0\le n < M_1$.
Hence $\sS_k =\emptyset$ for all $k \in \N$.
\item
$z_n \notin V_0 \cup\cV$ for $0 \le n < M_1$, 
and there exists a $k_0 \in \N$ such that 
$\sS_k:=[M_{k},M_{k+1})$ is an S-cycle for each $k \le k_0$ 
and $\sS_k=\emptyset $ for all $k >k_0$.   
\item
$z_n \notin V_0 \cup\cV$ for $0 \le n < M_1$, 
and $\sS_k:=[M_{k},M_{k+1})$ is a finite S-cycle for any $k \in \N$.
\end{enumerate}
Set $f=f_c$. For all cases, we have 
\begin{align*}
\sum_{n=1}^N\frac{1}{|Df^n(z)|}
&=
\sum_{n=1}^{M_1}\frac{1}{|Df^n(z)|}
+
\sum_{k \ge 1, \sS_k \neq \emptyset}
\sum_{i=1}^{M_{k+1}-M_k}
\frac{1}{|Df^i(z_{M_k})|~|Df^{M_k}(z)|}\\
&\le
\kappa_{\mathrm B}
+
\sum_{k \ge 1, \sS_k \neq \emptyset}
\frac{\kappa_{\mathrm A}}{|Df^{M_k}(z)|}
\end{align*}
by Lemmas A' and B'.
By Lemma B' again, we obviously have $|Df^{M_1}(z)|^{-1}<\kappa_{\mathrm B}$. 
Hence by Lemma C', we have 
$$
|Df^{M_k}(z)|
=
|Df^{M_k-M_{k-1}}(z_{M_{k-1}})|
\cdots
|Df^{M_2-M_{1}}(z_{M_1})|~|Df^{M_1}(z)|
\ge \lam^{k-1}/\kappa_{\mathrm B}.
$$
Hence we have
$$
\sum_{n=1}^N\frac{1}{|Df^n(z)|}
\le 
\kappa_{\mathrm B}
+
\sum_{k\ge 1}\frac{\kappa_{\mathrm A}\kappa_{\mathrm B}}{\lam^{k-1}}
<
\kappa_{\mathrm B}
+\kappa_{\mathrm A}\kappa_{\mathrm B}\frac{\lam}{\lam -1}
=:K_{\mathrm B}.
$$
\QED

\section{Hyperbolic metrics}\label{sec:Hyperbolic metrics}
For the proofs of Lemmas A', B', C' and C, we will use the hyperbolic metrics
and the expansion of $f_c$ with respect to these metrics.

For a domain $\Omega$ in $\C$ with $\#(\C-\Omega)\ge 2$,
there exists a {\it hyperbolic metric} $\rho(z)|dz|$ on $\Omega$ 
of constant curvature $-4$
induced by the metric $|dz|/(1-|z|^2)$ 
on the universal covering $\D = \widetilde{\Omega}$.
We first recall the following standard fact: 

\paragraph{Lemma W.}
{\it 
Let $\Omega_0$ be a domain in $\C$ with $\#(\C-\Omega_0)\ge 2$ and 
$\rho_0(z)|dz|$ be its hyperbolic metric.
Then for any domain $\Omega \subset \Omega_0$, 
the hyperbolic metric $\rho(z)|dz|$ of $\Omega$ satisfies
$$
\rho_0(z) \le \rho(z) \le \frac{1}{\dist(z,\partial \Omega)},
$$
where $\dist(z, \partial \Omega)$ is the Euclidean distance between 
$z$ and $\partial \Omega$. 
}

\medskip
See \cite[Theorems 1.10 \& 1.11]{Ah} for more details.

\paragraph{Postcritical sets.}
The {\it postcritical set} $P(f_c)$ of the polynomial $f_c(z)=z^2+c$ is
defined by
$$
P(f_c):=
\overline{
\brac{f_c(0), \,f_c^2(0),\, f_c^3(0),\, \cdots}
}.  
$$
For example, we have 
$$
P(f_\chat)=\{\hat{b}_1,\,\hat{b}_2,\, \cdots, \hat{b}_{l-1}\} \cup \Omega(\chat)
$$
when $c=\chat$ 
and this set is finite if $\chat$ is a Misiurewicz point.
Moreover, for any $c \approx \chat$, 
we have $\sharp P(f_c) \ge 2$ and 
the universal covering of (each component of) $\C-P(f_c)$
is the unit disk
\footnote{Without the parameter ray condition,
$f_c$ may have Siegel disks and the set $\C-P(f_c)$ 
may contain the disks.}. 

Let $\gamma=\gamma(z)|dz|$ denote the hyperbolic metric of $\C-P(f_\chat)$,
which is induced by the metric $|dz|/(1-|z|^2)$ on the unit disk $\D$. 
The metric $\gamma=\gamma(z)|dz|$ has the following properties: 
\begin{enumerate}[(i)]
\item
$\gamma: \C-P(f_\chat) \to \R_+$ is real analytic
and diverges on $P(f_\chat) \cup \{\infty\}$.
\item
if both $z$ and $f_\chat(z)$ are in $\C-P(f_\chat)$, 
we have 
$$
\frac{\gamma(f_\chat(z))}{\gamma(z)}|Df_\chat(z)| >1.
$$
\end{enumerate}

\paragraph{Lemma X.}
{\it If the constant $\nu$ is sufficiently small,
there exists a constant 
$C_{\mathrm X} \asymp \nu^2$ with the following property:
For any $c \approx \chat$, we have
$$
\frac{\gamma(z)}{\gamma(\zeta)} \ge C_{\mathrm X}
$$
if either 
\begin{enumerate}[\rm (1)]
\item
$z , \zeta \in J(f_c)-\cV$; or
\item
$z  \in J(f_c)-V_0 \cup \cV$ and $\zeta \in V_1-f_c(V_0)$.
\end{enumerate}
}

\paragraph{Proof.}
We may assume that there exists an $R_0>0$ such that 
$J(f_c)  \subset \overline{\D(R_0)}$ for any $c \approx \chat$.
Since $\gamma$ diverges only at the postcritical set $P(f_\chat)$ 
in $\overline{\D(R_0)}$,
there exists a constant $C_4>0$ such that $\gamma(w) \ge C_4$
for any $w \in \overline{\D(2 R_0)} - P(f_{\chat})$.
In particular, we have $\gamma(z) \ge C_4$ 
in both cases (1) and (2).
Moreover, for these cases,
we can find a constant $C_5$ independent of 
$\nu \ll 1$ and $c \approx \chat$ 
such that
$$
\dist (\zeta, P(f_\chat)) \ge C_5 \nu^2.
$$
Hence if $\nu$ is sufficiently small, then Lemma W implies that 
that $\gamma(\zeta) \le 1/(C_5\nu^2)$. 
Now we have $\gamma(z)/\gamma(\zeta) \ge C_4 C_5 \nu^2 =:C_{\mathrm X}$.\QED

\paragraph{Lemma Y.}
{\it There exists a constant $A>1$ such that 
for $c \approx \chat$, if $z, f_c(z), \ldots ,f_c^n(z)$ are all contained in 
$J(f_c)-\cV$, we have 
$$
|Df_c^n(z)| \ge C_{\mathrm X}A^n.
$$
This estimate also holds if $z, f_c(z), \ldots ,f_c^{n-1}(z)$ 
are all contained in $J(f_c)-V_0 \cup \cV$ and $f_c^n(z) \in V_1-f_c(V_0)$. 
}

\paragraph{Proof.}
Since the Julia set is uniformly bounded when $c \approx \chat$,
we may assume that there exists a constant $A>1$ such that
for any $c \approx \chat$,
$$
\frac{\gamma(f_c(w))}{\gamma(w)}|Df_c(w)| \ge A
$$
if either $w, \, f_c(w) \in J(f_c)-\cV$;
or $w  \in J(f_c)-\cV \cup V_0$ and $f_c(w) \in V_1-f_c(V_0)$.

By the chain rule, we have
$$
|Df_c^n(z)| 
= 
\prod_{i=0}^{n-1}|Df_c(f_c^i(z))|
\ge  
\prod_{i=0}^{n-1}\frac{\gamma(f_c^{i}(z))}{\gamma(f_c^{i+1}(z))}A
\ge
\frac{\gamma(z)}{\gamma(f^n_c(z))}A^n.
$$ 
By applying Lemma X with $\zeta:=f^n_c(z)$,
 we obtain the desired inequality. \QED

\section{Proof of Lemma B'}
Set $f=f_c$.
Suppose that $M<\infty$.
Since we have $z_i \notin V_0 \cup \cV$ for all $i \le M-1$, 
we can apply Lemma Y and we have
$$
|Df^i(z_0)| 
\ge \frac{\gamma(z_0)}{\gamma(z_i)} \cdot A^i
\ge C_{\mathrm X} A^i.
$$
If $z_{M} \notin  V_0 \cup \cV$ or $z_{M} \notin  V_1-f_c(V_0)$,
then we can apply Lemma Y again and we have 
$|Df^M(z_0)| \ge C_{\mathrm X} A^M \ge C_{\mathrm X}$.
Otherwise $z_{M} \in V_j$ for some $j \neq 1$.
Since $z_{M-1} \notin V_0 \cup \cV$,
we may assume that $|z_{M-1}| \ge \xi_0$ 
for some constant $0<\xi_0 \le 1/2$ depending only on $\chat$ and 
independent of $\nu \ll 1$,
$c \approx \chat$, and $z_0 \in J(f_c)$. 
Hence we have 
$$
|Df^M(z_0)|=|Df^{M-1}(z_{M-1})|~|Df(z_{M-1})| 
\ge C_{\mathrm X} A^{M-1}\cdot 2 \xi_0  \ge 2 \xi_0 C_{\mathrm X}.
$$
Thus 
$$
\sum_{i=1}^M\dfrac{1}{|Df^i(z_0)|}
\le 
\sum_{i=1}^{M-1}\dfrac{1}{C_{\mathrm X} A^i}
+
\dfrac{1}{2 \xi_0 C_{\mathrm X}}
<
\frac{1}{C_{\mathrm X}}
\paren{
\frac{1}{A-1} + 
\dfrac{1}{2 \xi_0}
}
=:\kappa_{\mathrm B}.
$$

If $M=\infty$, then the same estimate as above yields
$$
\sum_{i=1}^\infty\dfrac{1}{|Df^i(z_0)|}
\le 
\sum_{i=1}^ \infty \dfrac{1}{C_{\mathrm X} A^i}
<
\frac{1}{C_{\mathrm X} (A-1)}
<\kappa_{\mathrm B}.
$$

\QED

\section{Proof of Lemma A'}
Set $f=f_c$.
For a given S-cycle $\sS=[M,M')$,
we may assume that $M=0$ without loss of generality.
We divide the proof in two cases.

\paragraph{Case 1.}
Suppose that $\sS$ is either a finite S-cycle or an infinite S-cycle of type (I). 
Then there exist $j \in \{1,2, \cdots, l\}$, $m \in \N$,
and $L \in \N \cup \brac{\infty}$ such that 
\begin{itemize}
\item
$z=z_0 \in V_j$; 
\item
$z_{n-p} \in U_l$ when $n=(l-j)+mp$,  but $z_n \notin U_l$;
\item
$z_{n+i} \notin V_0 \cup \cV$ if $0 \le i <L$.   
\item
$M'<\infty$ iff $L<\infty$ and $M'=(l-j)+mp+L$.
\end{itemize}
Hence we have the following estimates of $|Df^n(z)|$:
\begin{itemize}
\item
When $n=1, \cdots, l-j-1$, we have
$z_n \in V_{j+n}$ and 
$$
|Df^n(z)| \ge \xi^n \ge \xi^{l-1}.  
$$
\item
When $n=(l-j) + kp+ i$ with $0 \le k < m$ and $0 \le i < p$,
\begin{align*}
|Df^n(z)| &= 
|Df^{l-j}(z)|~|Df^{kp}(z_{l-j})|~|Df^{i}(z_{(l-j)+kp})|\\
&\ge \xi^{l-j} \cdot \mu^k \cdot \xi^{i} \\
&\ge \xi^{(l-1)+(p-1)}\mu^k.
\end{align*}
\item
When $n=(l-j) + mp+ i$ with $0 \le i < L \le \infty$,
\begin{align*}
|Df^n(z)| &= 
|Df^{(l-j)+mp}(z)|~|Df^{i}(z_{(l-j)+mp})|\\
&\ge \xi^{l-j} \cdot \mu^m \cdot \frac{\gamma(z_{(l-j)+mp})}{\gamma(z_n)} 
\cdot A^i\\
&\ge \xi^{l-1}  \,C_{\mathrm X} \, A^i.
\end{align*}
Here the constant $A$ above is the same as that of Lemma Y.
\item
When $L<\infty$ and $n=M'=(l-j) + mp+ L$, 
the point $z_{M'}$ satisfies either
 $z_{M'} \in V_1-f_c(V_0)$; or
 $z_{M'} \in V_j$ for some $j \neq 1$.
By the same argument as in the proof of Lemma B', 
there exists a constant $0<\xi_0 \le 1/2$ depending only on $\chat$ such that
\begin{align}
|Df^{n} (z)|=|Df^{M'} (z)| &= 
|Df^{(l-j)+mp}(z)|~|Df^{L}(z_{(l-j)+mp})|\\
&\ge \xi^{l-j} \cdot 
\mu^m \cdot
\min
\{C_{\mathrm X} \, A^{L},\, C_{\mathrm X} \, A^{L-1} \cdot  2 \xi_0\}\notag \\
&\ge 2 \xi^{l-1} \,\xi_0 \, C_{\mathrm X}
\label{eq_B'}
\end{align}
\end{itemize}
By these estimates, when $M'<\infty$, we have:
\begin{align*}
&\sum_{i=1}^{M'} \frac{1}{|Df^i(z)|}\\
=&
\sum_{i=1}^{l-j-1} 
\frac{1}{|Df^i(z)|}
+
\sum_{k=0}^{m-1} 
\sum_{i=0}^{p-1} 
\frac{1}{|Df^{l-j}(z)|~|Df^{kp+i}(z_{l-j})|}\\
&~~~~~~~~~~~~~~~~~
+ 
\sum_{i=0}^{L-1} 
\frac{1}{|Df^{(l-j)+mp}(z)|~|Df^{i}(z_{(l-j)+mp})|}
+
\frac{1}{|Df^{M'}(z)|}
\\
\le& 
\frac{l-2}{\xi^{l-1}}
+
\sum_{k=0}^{m-1} 
\frac{p}{\xi^{(l-1)+(p-1)} \cdot \mu^k}
+ 
\sum_{i=0}^{L-1} 
\frac{1}{\xi^{l-1} \,C_{\mathrm X} \, A^i}
+
\frac{1}{ 2 \xi^{l-1} \,\xi_0 \, C_{\mathrm X}}
\\
\le&
\frac{l-2}{\xi^{l-1}}
+
\frac{p}{\xi^{(l-1)+(p-1)}}
\cdot \frac{\mu}{\mu-1}
+ 
\frac{1}{\xi^{l-1}\,C_{\mathrm X}}
\cdot 
\frac{A}{A-1}
+
\frac{1}{2 \xi^{l-1} \,\xi_0 \, C_{\mathrm X}}
\\
=:&\kappa_{\mathrm A}.
\end{align*}
Note that $\kappa_{\mathrm A}$
does not depend on $j,\,m$, and $L$.

If $M'=\infty$, then $L=\infty$ and one can easily check 
\begin{align*}
\sum_{i=1}^{\infty} \frac{1}{|Df^i(z)|}
\le
\frac{l-2}{\xi^{l-1}}
+
\frac{p}{\xi^{(l-1)+(p-1)}}
\cdot \frac{\mu}{\mu-1}
+ 
\frac{1}{\xi^{l-1} \, C_{\mathrm X}}
\cdot 
\frac{A}{A-1}<\kappa_{\mathrm A}.
\end{align*}

\paragraph{Case 2.}
Suppose that $\sS=[0,\infty)$ is an infinite S-cycle of type (II). 
Then there exists a $j \in \{1,2, \cdots, l\}$ such that $z=z_0 \in V_j$
and $z=b_j(c)$ if $j<l$ and $z \in \Omega(c)$ if $j=l$.
Hence for any $k \in \N$ we have $z_{(l-j)+kp} \in U_l$. 
By the same estimates as in Case 1, we have
\begin{align*}
\sum_{i=1}^{\infty} \frac{1}{|Df^i(z)|}
&=
\sum_{i=1}^{l-j-1} 
\frac{1}{|Df^i(z)|}
+
\sum_{k=0}^{\infty} 
\sum_{i=0}^{p-1} 
\frac{1}{|Df^{l-j}(z)|~|Df^{kp+i}(z_{l-j})|}\\
&
\le 
\frac{l-2}{\xi^{l-1}}
+
\sum_{k=0}^{\infty} 
\frac{p}{\xi^{(l-1)+(p-1)} \cdot \mu^k}\\
&=
\frac{l-2}{\xi^{l-1}}
+
\frac{p}{\xi^{(l-1)+(p-1)}} \cdot\frac{\mu}{\mu-1} <\kappa_A. 
\end{align*}
\QED

\section{Proof of Lemma C'}
Set $f=f_c$.
We will show that $|Df^{M'-M}(z_M)| \ge \kappa_{\mathrm{C}}/\nu$
for some constant $\kappa_{\mathrm{C}}$ that depends only on $\chat$. 
By choosing $\nu$ sufficiently small, 
we have $ \lam:=\kappa_{\mathrm{C}}/\nu>1$. 

As in the proof of Lemma A', we assume that $M=0$ and 
set $M':=(l-j)+mp +L$ where $z_0 \in V_j$ for some $1\le j\le l$.
We also set $n:=(l-j)+mp$, then by the chain rule we have
\begin{equation}
|Df^{M'}(z_0)|
=|Df^{n}(z_0)|\cdot |Df^{L}(z_n)|.
\label{eq_C'-0}
\end{equation}

First let us give an estimate of $|Df^{n}(z_0)|$.
We can find an $\tilde{R}_l>0$ such that 
$$
 f^p(\D(\hat{w},R_l)) \Subset \D(\hat{w}_p,\tilde{R}_l/2) \Subset \D(\hat{w}_p,\tilde{R}_l)
$$
for any $\hat{w}\in\Omega(\chat)$
if we choose $R_l$ small enough, where $\hat{w}_p=f_\chat^p(\hat{w})$. 
Let $\hat{x}:=\hat{b}_j$ if $z_0\in V_j$ and $j\not= l$,
 or $\hat{x}:=\hat{w}$ if $z_0\in\D(\hat{w},C_0'\nu^2) \subset V_l$ for some $\hat{w}\in\Omega(\chat)$. (The choice of $\hat{w}$ is not unique.) Let $x_0(c)=b_j(c)$ if $j<l$, or $x_0(c)=\chi_c(\hat{x})$ if $j=l$. 
Note that for any $c\approx \chat$, we have $b_j(c)\in V_j$, $\chi_c(\hat{w})\in \D(\hat{w},C_0'\nu^2)$, and $\chi_c(\hat{w}_p)\in\D(\hat{w}_p,C_0'\nu^2)$. In particular, we may assume that $|z_0-x_0(c)|\le \max(C_0,2C_0')\cdot \nu^2$ and $|\hat{x}_n-x_n(c)|\le R_l/2$, where $\hat{x}_n=f_\chat^n(\hat{x})$ and $x_n(c)=\chi_c(\hat{x}_n)=f_c^n(x_0(c))$. Thus, $|z_n-x_n(c)|\ge |z_n-\hat{x}_n|-|\hat{x}_n-x_n(c)|\ge R_l/2$.

Now we take the inverse branch $G$ of $f^n$
defined on $\D(\hat{x}_n,\tilde{R}_l)$ that maps $x_n(c)$ to $x_0(c)$, and $z_n$ to $z_0$. 
By the Koebe distortion theorem, we have
$$
|DG(z_n)| \asymp |DG(x_n(c))|
$$ 
and
$$ 
|z_0-x_0(c)|=|G(z_n)-G(x_n(c))| \asymp |DG(x_n(c))||z_n-x_n(c)|.  
$$
Since $|z_0-x_0(c)| \le \max(C_0\,\nu^2, 2C_0'\,\nu^2)$ and $|z_n-x_n(c)| \ge R_l/2$ ,
we have $|DG(z_n)| \le C_6\, \nu^2/R_l$, where $C_6$ is a constant 
independent of $c \approx \chat$, $\nu \ll 1$, and $z_0 \in J(f_c)$. Hence $|Df^n(z_0)| \ge R_l/(C_6\,\nu^2)$.

Next we give an estimate of the form $|Df^{L}(z_n)| \ge C_7 \nu$,
 where $C_7$ is a constant 
independent of $c \approx \chat$, $\nu \ll 1$, and $z_0 \in J(f_c)$.
(Then by (\ref{eq_C'-0}) the proof is done.)
The estimate relies on the geometry of (and dynamics on) the postcritical set $P(f_\chat)$:
Take any $i \in [0, L)$, then by Lemmas W and Y we obtain 
\begin{align*}\label{eq_C'-0}
|Df^{L}(z_n)|
=&|Df^{L-i}(z_n)||Df^i(z_{M'-i})|\\
\ge& 
\frac{\gamma(z_n)}{\gamma(z_{M'-i})} A^{L-i}
\cdot |Df^i(z_{M'-i})|\\
\ge& 
\gamma(z_n)\cdot 
\dist(z_{M'-i}, P(f_\chat))
\cdot |Df^i(z_{M'-i})|.
\end{align*}
By taking a small enough $R_l$, we may assume that 
$f_c^p(U_l)$ is disjoint from $P(f_\chat)-\Omega(\chat)$.
Hence $z_n$ has a definite distance from $P(f_\chat)$
(more precisely, $\dist(z_{n}, P(f_\chat))$ is bigger than 
a positive constant independent of $c \approx \chat$,
$\nu \ll 1$, and $z_0 \in J(f_c)$)
and we always have $\gamma(z_n) \asymp 1$.

Thus it is enough to show: 
{\it There exists an $i \in [0,\,l+p)$ such that
\begin{itemize}
\item[\rm (1)]
$z_{M'-i}$ has a definite distance from $P(f_\chat)$; and
\item[\rm (2)]
$|Df^i(z_{M'-i})|\ge C_8 \nu$ for some constant $C_8$ depending only on $\chat$. 
\end{itemize}
} 
Note that if $z_{M'} \in V_{0}$, then $z_{M'}$ already has a definite distance from $P(f_\chat)$ by semi-hyperbolicity. This situation corresponds to $i=0$ and condition (2) is ignored.

\begin{figure}[htbp]
\begin{center}
\includegraphics[width=.75\textwidth]{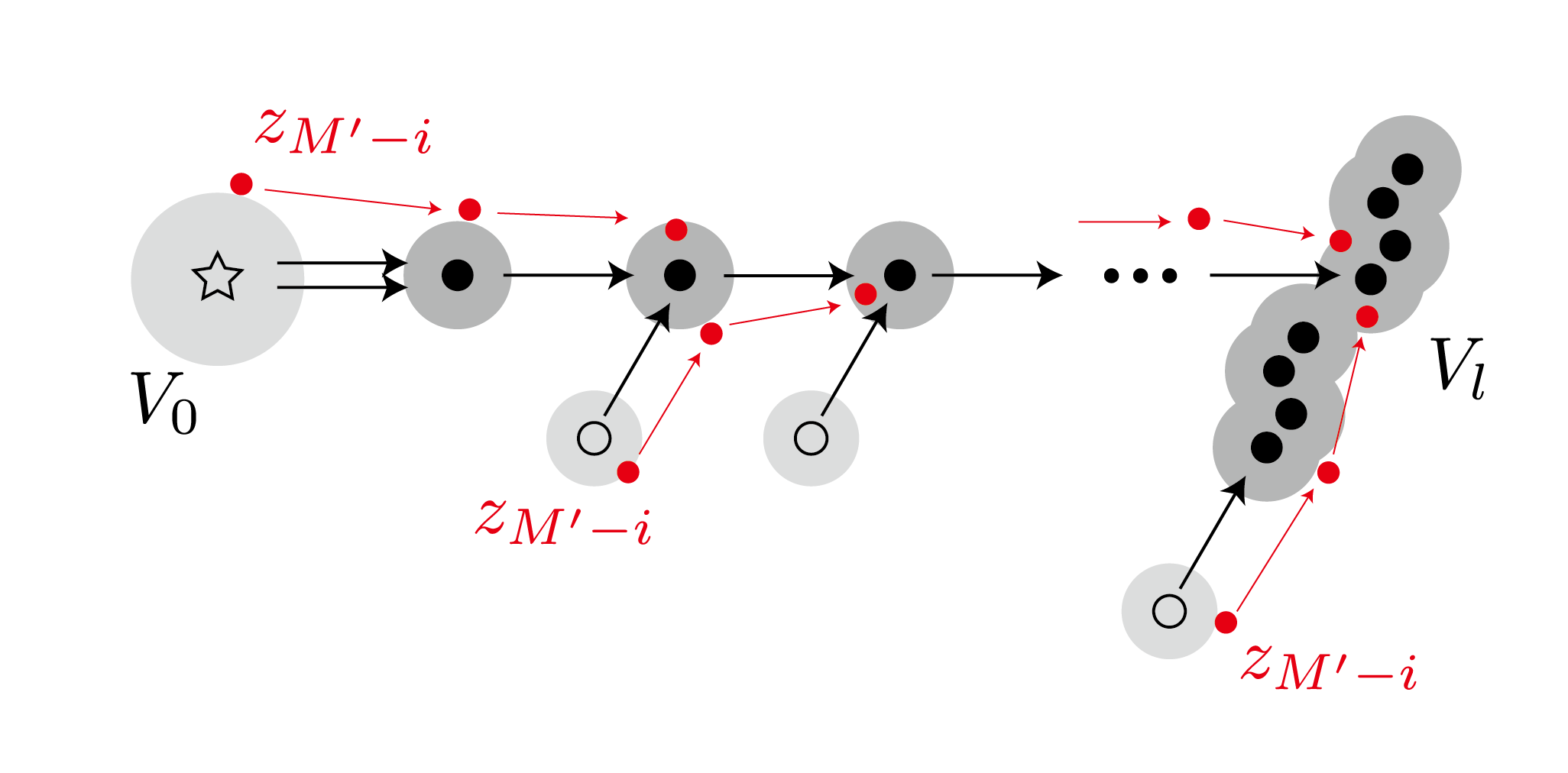}
\caption{Black heavy dots indicate the critical orbit. 
Some possible behaviors of $z_{M'-i} \mapsto z_{M'-i+1} \mapsto \cdots \mapsto z_{M'}$ 
are indicated by smaller dots (in red). 
}
\label{fig_Postcritical}
\end{center}
\end{figure}

If $z_{M'} \in V_{j'}$ with $1 \le j' \le l$,
then such an $i$ can be found in $[1,l+p)$ by the following procedure (Figure \ref{fig_Postcritical}). 
Suppose that $z_{M'} \in V_{1}$. Then $z_{M'-1}$ is contained in $f^{-1}(V_1)-V_0$, and thus $|z_{M'-1}| \ge \nu$. By setting $i=1$, it follows that $z_{M'-1}$ has a definite distance from $P(f_\chat)$, and we have $|Df(z_{M'-1})| \ge 2 \nu$.

Suppose that $z_{M'} \in V_{2}$. Then $f^{-1}(V_{2})$ has two components
containing $\pm \hat{b}_{1}$ for any $c \approx \chat$.
If $z_{M'-1}$ is in the component containing $-\hat{b}_{1}$,
then $|z_{M'-1}-(- \hat{b}_{1})| \asymp \nu^2$ 
and it has a definite distance from $P(f_\chat)$.
Now set $i=1$.
Since 
$
|Df(-\hat{b}_{1})|= 2 |\hat{b}_{1}|\ge \xi 
$
by definition of $\xi$ in Section 2, we have 
$|Df(z_{M'-1})| \asymp |Df(-\hat{b}_{1})| \ge  \xi>\nu$ 
for $\nu \ll 1$.
If $z_{M'-1}$ is in the component containing $\hat{b}_{1}$, then $z_{M'-1}$ is necessarily contained in $f^{-1}(V_2)-V_1$, and
then $|z_{M'-1}-\hat{b}_{1}| \asymp \nu^2$.
In this situation $|z_{M'-2}| \asymp \nu$ and $z_{M'-2}$ 
has a definite distance from $P(f_\chat)$. 
Set $i=2$. Then 
$$
|Df^2(z_{M'-2})|=|Df(z_{M'-2})||Df(z_{M'-1})| 
\asymp \nu \cdot |Df(\hat{b}_{1})| \ge \xi  \nu .
$$

Suppose that $z_{M'} \in V_{j'}$ with $j'=3, \, \cdots, l-1$.
As in the situation of $z_{M'} \in V_{2}$, either
\begin{itemize}
\item
$|z_{M'-i} - (-\hat{b}_{j'-i})| \asymp \nu^2$ 
for some $i < j'$ and $z_{M'-i}$ has a definite distance from $P(f_\chat)$; or
\item
$|z_{M'-j'}| \asymp \nu$ and $z_{M'-j'}$ has a definite distance from $P(f_\chat)$. We set $i:=j'$ in this case.
\end{itemize}
In both cases, we have $|z_{M'-k} - \hat{b}_{j'-k}| \asymp \nu^2$ for each $k=1, \cdots, i-1$.
In particular, since $2|\hat{b}_{n}| \ge \xi$ for $n \in \N$, 
we have:
\begin{itemize}
\item
If $i<j'$, then 
$
|Df^i(z_{M'-i})| \asymp
2 |\hat{b}_{j'-i}| \cdot 2 |\hat{b}_{j'-i+1}| \cdots  2 |\hat{b}_{j'-1}|  
\ge \xi^i \ge \xi^{l-2}. 
$
\item
If $i=j'$, then 
$
|Df^i(z_{M'-i})| \asymp
 2 \nu \cdot 2 |\hat{b}_{1}| \cdots 2 |\hat{b}_{j'-1}| 
\ge 2 \xi^{j'-1} \nu  \ge 2 \xi^{l-2} \nu. 
$
\end{itemize}
In both cases, we have $|Df^i(z_{M'-i})| \ge C_8 \nu$
for some constant $C_8>0$ independent of $c \approx \chat$,
$\nu \ll 1$, and $z_0 \in J(f_c)$.

Finally suppose that $z_{M'} \in V_{l}$, i.e., 
$\dist (z_{M'}, \Omega(\chat)) < C_0'\, \nu^2$ by definition of $V_l$. 
Now we claim: 
{\it 
there exists a  $k'\le p$ such that $\dist (z_{M'-k'}, \Omega(\chat))\asymp R_l$.}

Indeed, if there exists some $1 \le k'<p$ such that 
$z_{M'},\, z_{M'-1},\, \cdots,\,  z_{M'-k'+1} \in U_l$ but $z_{M'-k'}\not\in U_l$, 
then $\dist (z_{M'-k'}, \Omega(\chat))\asymp R_l$. 
Now suppose that 
all $z_{M'}$, $z_{M'-1}$, $z_{M'-2}, \ldots, z_{M'-p+1}$ remain in $U_l$ 
(but not in $V_l$ except $z_{M'}$).
Let us show that $z_{M'-p} \notin U_l$ by contradiction.

Assume that $z_{M'-p} \in U_l$.
Since $|Df^p(z)| \ge \mu > 2.5$ for $z \in U_l$, 
by the Koebe distortion theorem and invariance of $\Omega(c)$ by $f^p=f_c^p$,
we obtain $2\cdot \dist (z_{M'-p}, \Omega(c)) < \dist (z_{M'}, \Omega(c))$ if $\nu \ll 1$.
(Note that we have $z_{M'} \notin \Omega(c)$ since $\sS$ is a finite S-cycle.)
Since $\Omega(\chat)$ moves holomorphically,
we may assume that
$\dist( \Omega(c), \Omega(\chat)) \le C_0'\nu^2/4$ for $c \approx \chat$.
Hence we obtain
\begin{align*}
\dist (z_{M'-p}, \Omega(\chat)) 
&\le \dist (z_{M'-p}, \Omega(c)) +C_0'\nu^2/4\\
&< \dist (z_{M'}, \Omega(c))/2  + C_0'\nu^2/4 \\
&\le (\dist (z_{M'}, \Omega(\chat)) +C_0'\nu^2/4)/2 + C_0'\nu^2/4\\
&< C_0'\nu^2.
\end{align*}
It would imply $z_{M'-p}\in V_l$, contradicting the construction of the S-cycle $[M, M')$.
It follows that $z_{M'-p}\not\in U_l$ and
thus $\dist (z_{M'-k'}, \Omega(\chat))\asymp R_l$ for $k'=p$.

The point $z_{M'-k'}$ above has a definite distance from $\Omega(\chat)$. 
It also has a definite distance from $P(f_\chat)$,
unless 
$|z_{M'-k'} - \hat{b}_{l-1}| \asymp \nu^2$.
However, in this case we may apply the same argument as 
in the case of $1 \le j' \le l-1$ and 
there exists an $i \in [k', k'+l)$ such that 
$z_{M'-i}$ has a definite distance from $P(f_\chat)$. 
Moreover, since $i$ is bounded by $p+l$, 
we have $|Df^i(z_{M'-i})| \ge C_8 \nu$ 
by replacing the above $C_8$ if necessary.  
\QED 

\section{Proof of Lemma C}
This proof is similar to that of Lemma C'.
We will show that $|Df_c^{N'-N}(z_0)| \ge K_{\mathrm{C}}/\nu$
for some constant $K_{\mathrm{C}}$ that depends only on $\chat$,
and we set $ \Lam:=K_{\mathrm{C}}/\nu>1$ by choosing $\nu \ll 1$.  

Without loss of generality we may assume that $N=0$.
Set $n:=l+mp$ and $L:=N'-n$
such that $z_{n-p} \in U_l$, $z_n \notin U_l$,
$z_{n+i} \notin V_0 $ for $0 \le i < L$, 
and $z_{n+L} \in V_0$.

By the chain rule, we have 
\begin{equation}
|Df_c^{N'}(z_0)|
=|Df_c^{n}(z_0)|\cdot |Df_c^{L}(z_n)|.
\label{eq_C-0}
\end{equation}
By Lemma U, we have 
$|Df_c^n(z_0)| 
\ge C_{\mathrm U}/|z_0| 
\ge C_{\mathrm U}/\nu
$ 
where the constant $C_{\mathrm U}>0$ is 
independent of $c \approx \chat$ and $z_0 \in J(f_c) \cap V_0$. 
Hence it is enough to show that $|Df_c^{L}(z_n)| \ge \eta$
for some constant $\eta>0$ that is independent of $\nu \ll 1$, 
$c \approx \chat$ and $z_0 \in V_0 \cap J(f_c)$. 
(Then we have $|Df_c^{N'}(z)| \ge C_{\mathrm U}\eta/\nu$ by (\ref{eq_C-0}) 
and the proof is done.)

To show this, we use the hyperbolic metric.
Let $\rho(z)|dz|=\rho_c(z)|dz|$ be the hyperbolic metric on $\C-P(f_c)$,
where 
$$
P(f_c)=\brac{c,\,f_c(c),\,f_c^2(c),\,\ldots} 
$$ 
is the postcritical set of $f_c$ for $c \approx \chat$ with $c \notin \M$.

Since $J(f_c) \cap P(f_c) =\emptyset$ when $c \notin \M$,
we have 
$$
\frac{\rho(f(z))}{\rho(z)}|Df_c(z)| \ge 1
$$
for any $z \in J(f_c)$. (See \cite[Theorem 3.5]{Mc} for example.)
We also have $\rho(z) \le \dist(z, P(f_c))^{-1}$ by Lemma W.
Hence $|Df^L(z_n)|\ge \rho(z_n)/{\rho(z_{N'})} 
\ge \rho(z_n)\cdot\dist(z_{N'},P(f_c))$. 
To complete the proof, 
we show that both $\rho(z_n)$ and $\dist(z_{N'},P(f_c))$ are 
uniformly bounded from below 
for any $c \approx \chat$ and for any $z_0 \in V_0=\D(\nu)$.

Let us work with $\dist(z_{N'},P(f_c))$ first:
Let $\widetilde{\cR}(c)$ denote
the closure of the union of the 
forward images of the dynamic ray $\cR_c(\theta)$.
By using the set $\widehat{\cR}(c)$ defined in Section \ref{sec:ParameterRayCondi},
we have
$$
\widetilde{\cR}(c) =\overline{
\cR_c(\theta) \cup \cR_c(2\theta) \cup \cdots \cup \cR_c(2^{l-1} \theta) 
\cup \widehat{\cR}(c)}.
$$
By 
Proposition S, 
this set moves continuously 
as $c \to \chat$ along $c \in\cR_\M(\theta)$
with respect to the Hausdorff distance on the sphere. 
Since the postcritical set $P(f_c)$
is contained in $\widetilde{\cR}_c$, we obtain
$$
\dist(z_{N'},P(f_c))
\ge \dist(z_{N'},\widetilde{\cR}_c)
\ge \dist(0,\widetilde{\cR}_c)
-|z_{N'}|
\ge \dist(0,\widetilde{\cR}_c)-\nu,
$$
where $\dist(0,\widetilde{\cR}_c)$ tends to $\dist(0,\widetilde{\cR}_\chat)>0$
as $c \to \chat$ with $c \in \cR_{\M}(\theta)$.
Now we choose sufficiently small $\nu$ and 
we conclude that 
$\dist(z_{N'},P(f_c))$
is bounded by a positive constant that is 
independent of $c \to \chat$ with the parameter ray condition
and $z_{N'} \in V_0$.

Next we work with $\rho(z_n)$:
Let $T_c:\C \to \C~(c \neq 0)$ be a complex 
affine map with $T_c(c)=\chat$ and $T_c(f_c(c))=f_\chat(\chat)$
such that $T_c(z) \to z$ uniformly on compact sets as $c \to \chat$. 
Set $g_c:=T_c \circ f_c \circ T_c^{-1}$. 
Then $g_c$ is a quadratic map whose postcritical set is 
$$
P(g_c)=T_c(P(f_c))=
\brac{\chat,\,f_\chat(\chat)=g_c(\chat),\,g_c^2(\chat),\,\ldots}. 
$$ 
Hence the hyperbolic metrics $\rho_c'$ on $\C-P(g_c)$
and $\hat{\rho}$ on $\C-\brac{\chat,\,f_\chat(\chat)}$ satisfy 
$T_c^\ast \rho'_c=\rho_c$ and $\hat{\rho} \le \rho_c'$ for all $c$, where $T_c^\ast$ is the pull-back.

As in the proof of Lemma C',
if we choose $R_l$ small enough,  
then we can find an $\tilde{R}_l>0$ such that 
$
f_c^p(U_l) \Subset \mathrm{N}(\Omega(\chat),\tilde{R}_l)
$
for any $c \approx \chat$ and that 
the closure $E$ of the set 
$\mathrm{N}(\Omega(\chat),\tilde{R}_l)-V_l$ 
contains neither $\chat$ nor $f_\chat(\chat)$.
(Note that $f_\chat(\chat)$ may belong to $\Omega(\chat)$ and be contained in $V_l$.)
It follows that $z_n$ is contained in $E$ for $c \approx \chat$, 
and hence so is $z_n':=T_c(z_n)$. Thus we obtain 
$$
\rho'_c(z_n')\ge \hat{\rho}(z_n') \ge \min_{w \in E} \hat{\rho}(w) >0.
$$
Since $\rho_c(z_n)=\rho'(T_c(z_n))|DT_c(z_n)|=\rho'(z_n')|DT_c(z_n)|$ and 
$DT_c(w) \to 1$ uniformly on $E$ as $c \to \chat$,
we conclude that $\rho(z_n)$ is bounded by a positive 
constant from below that is independent of $\nu \ll 1$, $c \approx \chat$ 
and the original choice of $z_0 \in V_0 \cap J(f_c)$.
\QED

\section{Itinerary sequences} \label{sec:itinerary}

 When $c\not\in \M$, the critical value $c$  has a well defined external 
 angle $t_c=(2 \pi)^{-1}\arg \Phi_c(c)$.
The angle $t_c$ is not equal to zero when $c\in\mathbb{X}=\C-\M\cup\R_+$.  
For $c\in\mathbb{X}$, the dynamic rays $\mathcal{R}_c(t_c/2)$ 
and $\mathcal{R}_c((t_c +1)/2)$ together with the critical point $0$ separate the complex plane $\mathbb{C}$ into two disjoint open sets, say $W_0=W_0(c)$  and $W_1=W_1(c)$. Let the one that contains $c$ be $W_0$.
If $t_c=\theta$ and $\mathcal{R}_\M(\theta)$ lands at a semi-hyperbolic parameter $\chat$, then  $\mathcal{R}_\chat(\theta)$ lands at $\chat$, and both $\mathcal{R}_\chat(\theta/2)$ and $\mathcal{R}_\chat((\theta+1)/2)$ land at $0$. Moreover, as $c$ approaches $\chat$ along $\mathcal{R}_\M(\theta)$, in a large disk centered at $0$, rays $\mathcal{R}_c(\theta/2)$ and $\mathcal{R}_c((\theta+1)/2)$ move continuously to 
$\mathcal{R}_\chat(\theta/2)$ and $\mathcal{R}_\chat((\theta+1)/2)$, respectively.

Assume $z\in J(f_c)$. Define its {\it itinerary} or {\it itinerary sequence} $I_c(z)=\{I_c(z)_n\}_{n\ge 0}$ by $I_c(z)_n=0$ if $f_c^n(z)\in W_0$, $I_c(z)_n=1$ if $f_c^n(z)\in W_1$, and $I_c(z)_n=*$ if $f_c^n(z)=0$. If the critical point $0$ belongs to the Julia set, $I_c(f_c(0))$ is called the {\it kneading sequence} for $f_c$.  

\begin{remark}
 One can also define an itinerary $\{s_0, s_1,\ldots\}$ in such a way that $s_n=0$ if $f_c^n(z)\in \overline{W_0}\cap J(f_c)$ and that $s_n=1$ if $f_c^n(z)\in \overline{W_1}\cap J(f_c)$. When $\chat$ is a semi-hyperbolic parameter, if $f_\chat^k(z)=0$ for some $k\ge 0$, then $f_\chat^n(z)\not=0$ for all $n\not=k$ since the critical point  is non-recurrent. Suppose $f_\chat^n(z)\in W_{s_n}$ for $n\not= k$, then with the definition of itinerary in this remark, the itinerary of $z$ will have two values $\{s_0,\ldots,s_{k-1},0,s_{k+1},\ldots\}$ and $ \{s_0,\ldots,s_{k-1},1,s_{k+1},\ldots\}$. We employ the symbol $*$ in the above definition  so as to identify sequences $\{s_0,\ldots,s_{k-1},0,s_{k+1},\ldots\}$ and $\{s_0,\ldots,s_{k-1},1,s_{k+1},\ldots\}$ by the one $\{s_0,\ldots,s_{k-1},*,s_{k+1},\ldots\}$.
\end{remark}

\paragraph{Lemma Z.}
 Let $\chat$ be a semi-hyperbolic parameter. 
\begin{enumerate}[\rm (i)]
\item
$I_\chat(z)=I_\chat(w)$ if and only if $z=w$. 
\item
If $I_\chat(z)_k=*$ and $I_\chat(z)_n=I_\chat(w)_n$ for all $n\not=k$, 
then $I_\chat(w)_k=*$ and $w=z$.
\end{enumerate}
\paragraph{Proof.}
 Since $\mathbb{C}-\mathcal{R}_\chat(\theta)\cup\{\chat\}$ is a simply connected domain without a critical value, there exist inverse branches $f_{\chat, i}^{-1}:\mathbb{C}-\mathcal{R}_\chat(\theta)\cup\{\chat\} \to W_i$ of $f_\chat$, $i=0$ or $1$. Each of these two branches can be extended at the critical value $\chat$, and each extended branch is one-to-one.
 
(i) If $I_\chat(z)_n=I_\chat(w)_n=s_n$ for all $n\ge 0$, then for any $N\in\N_0$ both $f_\chat^N(z)$ and  $f_\chat^N(w)$ belong to $W_{s_N}$ provided $s_N\not=*$, or belong to $\{0\}$   provided $s_N=*$. The set $J(f_\chat)\cap \overline{W_{s_N}}$ can be convered by a finite number of disks $\D(y_i, \epsilon)$ with $y_i\in  J(f_\chat)\cap \overline{W_{s_N}}$, $i\in F$, and $F$ is a finite index set. We choose $\epsilon$ to be the constant such that the inequality \eqref{CJYcontraction} holds. Let $B_N(y_i, \epsilon)$ be the component of $f_\chat^{-N}(\D(y_i,\epsilon))$ such that $f_\chat^{-N}(y_i)\in \overline{W_{s_0}}$,  $f_\chat^{-N+1}(y_i)\in \overline{W_{s_1}}, ~\ldots~, ~ f_\chat^{-1}(y_i)\in \overline{W_{s_{N-1}}}$. It is not difficult to see that both $z$ and $w$ are contained in a simply connected domain covered  by the union $\bigcup_{i\in F}B_N(y_i,\epsilon)$. It follows that  $z=w$ easily from the exponential contraction \eqref{CJYcontraction} by taking $N\to\infty$. 

(ii) If $I_\chat(z)_k=*$ and $I_\chat(z)_n=I_\chat(w)_n$ for all $n>k$, then $I_\chat(f_\chat^{k+1}(z))=I_\chat(f_\chat^{k+1}(w))$. Thus, $f_\chat^{k+1}(z)=f_\chat^{k+1}(w)=\chat$ by (i). Since $\chat$ is the critical value, $f_\chat^k(w)=0$, and then $I_\chat(w)_k=*$. Therefore, $I_\chat(z)=I_\chat(w)$, and then $z=w$ by (i). 
\QED \\

Let $z(c)$ and $\chat$ be as in Theorem \ref{thm_main}, and let $c_0$ be $c(2)$ in 
Theorem \ref{thm_HolomorphicMotionLands} or be as in Theorem \ref{thm_C2SH}. The statement (i) of following corollary describes how the itinerary of $z(c)$ retains. The statement (ii) tells that every given point, say $w$, of $J(f_\chat)$ is a limiting point $z(\chat)$ of some $z(c)$ in $J(f_c)$ where the limit is taken as in Theorem \ref{thm_HolomorphicMotionLands}.  

\begin{cor}\label{corformain}
\hfill{}
\begin{itemize}
\item[\rm (i)] Suppose $I_{c_0}(z(c_0))={\bf s}$, then $I_\chat(z(\chat))={\bf s}$ if and only if $f_\chat^n (z(\chat))\not=0$ for all $n\ge 0$, otherwise $I_\chat(z(\chat))=\{s_0,\ldots,s_{k-1},*,s_{k+1},\ldots\}$ if and only if $f_\chat^k(z(\chat))=0$ for some $k\ge 0$. 
\item[\rm (ii)] Let $w \in J(f_\chat)$ and $I_\chat(w)={\bf s}$. If $f^n_\chat(w)\not=0$ for all $n\ge 0$, there exists a unique $z(c_0)$, with $I_{c_0}(z(c_0))={\bf s}$, such that $w = z(\chat)$. If $f^k_\chat(w)=0$ for some $k\ge 0$, then there exist exactly two $z(c_0)$ and $\tilde{z}(c_0)$, having itineraries $\{s_0, \ldots, s_{k-1},0,s_{k+1},\ldots\}$ and $\{s_0, \ldots, s_{k-1},1,s_{k+1},\ldots\}$ respectively, such that $w = z(\chat)=\tilde{z}(\chat)$.
\end{itemize}
\end{cor}
\paragraph{Proof.}
(i) For $c\not\in \M$, every point $z\in J(f_c)$ of given itinerary is bounded away from $\mathcal{R}_c(\theta/2)\cup\mathcal{R}_c((\theta+1)/2)\cup\{0\}$ and moves holomorphically with $c$. Thus, $f_\chat^n (z(\chat))\in \overline{W_{s_n}(\chat)}$ if $f_{c_0}^n(z(c_0))\in W_{s_n}(c_0)$. Hence, $I_\chat(z(\chat))={\bf s}$ if $f_\chat^n(z(\chat))\not=0$ for all $n\ge 0$. If $f_\chat^k(z(\chat))=0=\overline{W_0(\chat)}\cap\overline{W_1(\chat)}\cap J(f_\chat)$, then $0\not= f_\chat^n(z(\chat))\in W_{s_n}(\chat)$ for all $n\not=k$ and $I_\chat(z(\chat))=\{s_0,\ldots,s_{k-1},*,s_{n+1},\ldots\}$. 

(ii) For any $w\in J(f_\chat)$, by Theorem \ref{thm_C2SH}, there exists $z(c_0)\in J(f_{c_0})$ such that $h_c(z(c_0))=z(c)\to z(\chat)=w$ as $c\to \chat$ along $\mathcal{R}_{\mathbb{M}}(\theta)$. If $f_\chat^n(w)\not=0$ for all $n\ge 0$, then $z(c_0)\not=0$ for all $n\ge 0$, and we conclude that $I_{c_0}(z(c_0))=I_\chat(w)$. If there exists another $\tilde{z}(c_0)\in J(f_{c_0})$ such that $h_c(\tilde{z}(c_0))=\tilde{z}(c)\to \tilde{z}(\chat)=w$ as $c\to \chat$ along $\mathcal{R}_{\mathbb{M}}(\theta)$, then $I_{c_0}(\tilde{z}(c_0))=I_{c_0}(z(c_0))$, and consequently $\tilde{z}(c_0)=z(c_0)$ by the bijectivity between the itinerary sequences and Julia set $J(f_{c_0})$. 

If
$f_\chat^k(w)=0$ for some $k\ge 0$, then $f_{c_0}^n(z(c_0))\in W_{s_n}(c_0)$ for $n\not=k$, and $f_{c_0}^k(z(c_0))$ belongs to $W_0(c_0)$ or $W_1(c_0)$. Without loss of generality, assume $f_{c_0}^k(z(c_0))\in W_0(c_0)$. Let $\tilde{z}(c_0)$ be such a point that $f_{c_0}^{k+1}(\tilde{z}(c_0))=f_{c_0}^{k+1}(z(c_0))$, $f_{c_0}^k(\tilde{z}(c_0))\in W_1(c_0)$, and $f_{c_0}^n(\tilde{z}(c_0))\in W_{s_n}(c_0)$ for $0\le n<k$. It is easy to see that such a point exists. We have $f_{c_0}^{k+1}(\tilde{z}(c))\to f_{\chat}^{k+1}(w)$ as $c\to \chat$. And, by (i) and Lemma Z, we obtain $I_\chat(\tilde{z}(\chat))= I_\chat(w)$ and $\tilde{z}(\chat)=w$. If there is another $z'(c_0)\in J(f_{c_0})$ such that $h_c(z'(c_0))\to w$ as $c\to\chat$ along $\mathcal{R}_{\mathbb{M}}(\theta)$, then either $I_{c_0}(z'(c_0))=I_{c_0}(z(c_0))$ or $I_{c_0}(z'(c_0))=I_{c_0}(\tilde{z}(c_0))$. Consequently, by the bijectivity between the itinerary sequences and Julia set $J(f_{c_0})$, we conclude that $z'(c_0)=z(c_0)$ or $\tilde{z}(c_0)$.
\QED \\

\section{Proofs of Theorems \ref{corformain2} and \ref{thm_C2SH_2}}
\label{sec:Proofs of Theorems 1.5 and 1.6}

\paragraph{Proof of Theorem \ref{corformain2}.}
Because $J(f_\chat)$ is locally connected, it is clear that $\mathcal{E}^\theta (\theta)=I_\chat(\chat)$, namely the kneading sequence of $\theta$ is equal to the kneading sequence for $f_\chat$. Hence, it is enough to prove the theorem by using ${\bf e}=I_\chat (\chat)$. Note  that ${\bf e}\in \Sigma_2$ because $\chat$ is not recurrent under iteration of $f_\chat$.

 For any $w\in J(f_\chat)$, we have  $\sigma^n(I_\chat(w))\not={\bf e}$ for all $n\ge 0$, or $I_\chat(w)={\bf e}$, or $\sigma^k(I_\chat(w))={\bf e}$ for some $k\ge 1$. For any ${\bf s}\in\Sigma_2$ satisfying $\sigma^n({\bf s})\not={\bf e}$ for all $n\ge 0$ or ${\bf s}={\bf e}$, from Corollary \ref{corformain}, there corresponds a unique $w\in J(f_\chat)$ with $I_\chat(w)={\bf s}$. For such ${\bf s}\in\Sigma_2$ that    $\sigma^{k+1}({\bf s})={\bf e}$ for some $k\ge 0$, there is a unique ${\bf a}\not={\bf s}$ in $\Sigma_2$ satisfying ${\bf a}\sim_{\bf e}{\bf s}$ and again from Corollary \ref{corformain} there corresponds a unique $w\in J(f_\chat)$ with $I_\chat(w)=\{a_0,\ldots,a_{k-1},*,a_{k+1},\ldots\}= \{s_0,\ldots,s_{k-1},*,s_{k+1},\ldots\}$. This shows the bijectivity between $\Sigma_2/{\sim_{\bf e}}$ and $J(f_\chat)$. Let the bijection  $\Sigma_2/{\sim_{\bf e}} \to J(f_\chat)$ be $h$. Since $I_\chat(h({\bf s}))={\bf s}$ if $f_\chat^n( h({\bf s}))\not=0$ for all $n\ge 0$ or $I_\chat(h({\bf s}))=\{s_0,\ldots,s_{k-1},*,s_{k+1},\ldots\}$ if $f_\chat^k( h({\bf s}))=0$ for some $k\ge 0$ (we use ${\bf s}$ for an element in both $\Sigma_2$ and $\Sigma_2/{\sim_{\bf e}}$ if it does not cause any confusion), by a similar argument to the proof of Lemma Z (i), the continuity of $h$ follows easily by  virtue of the exponential contraction \eqref{CJYcontraction}. Compactness of  $\Sigma_2/{\sim_{\bf e}}$ and $J(f_\chat)$ leads to $h$ a homeomorphism. To show $h$ acts as a conjugacy, observe from Corollary \ref{corformain} that points $h\circ\sigma({\bf s})$ and $f_\chat\circ h({\bf s})$ have the same itinerary under $f_\chat$, thus they are the same by Lemma Z (i).
\QED

\paragraph{Proof of Theorem \ref{thm_C2SH_2}.}
There are exactly two cases: $f^n_{\chat}(w)\not=0$ for all $n\ge 0$ or $f^n_\chat(w)=0$ for some $n\ge 0$. By Corollary \ref{corformain}, $h^{-1}_\chat(\{w\})$ is a singleton if and only if $f^n_\chat(w)$ is as the first case, whereas it consists of two distinct points if and only if $f^n_\chat(w)$ is as the second case.
\QED

%%%%%%%%

\section*{Acknowledgments}
The authors would like to thank the referee for 
his/her comments and suggestions 
that make the paper more precise and readable. 
Chen was partly supported by NSC 99-2115-M-001-007, MOST 103-2115-M-001-009, 104-2115-M-001-007, and 105-2115-M-001-003. 
Kawahira was partly supported by JSPS KAKENHI Grant Number 16K05193. 
They thank the hospitality of Academia Sinica, Nagoya University, 
RIMS in Kyoto University, and Tokyo Institute of Technology
where parts of this research were carried out.

\vspace{1cm}

{~}\\
Yi-Chiuan Chen \\
Institute of Mathematics\\
Academia Sinica\\
Taipei 10617, Taiwan\\
YCChen@math.sinica.edu.tw

\vspace{.5cm}

{~}\\
Tomoki Kawahira\\
Department of Mathematics\\ 
Tokyo Institute of Technology\\
Tokyo 152-8551, Japan\\
kawahira@math.titech.ac.jp  \\

\noindent
Mathematical Science Team\\ 
RIKEN Center for Advanced Intelligence Project (AIP)\\
1-4-1 Nihonbashi, Chuo-ku\\ 
Tokyo 103-0027, Japan

\end{document}